\documentclass[preprint,12pt]{elsarticle}

\usepackage{amssymb}
\usepackage{amsthm}
\usepackage{graphics}
\usepackage{epsfig}
\usepackage{amssymb}
\usepackage{hyperref}
\usepackage{graphicx}
\usepackage{subfigure}
\usepackage{color}
\usepackage{tikz}
\usetikzlibrary{patterns}
\usepackage[ruled,vlined]{algorithm2e}
\usepackage{hyperref}
\usepackage{a4wide}
\usepackage{amsmath}
\usepackage{amsfonts}
\usepackage{enumerate}
\newcommand{\pf}{\noindent {\bf Proof: }}

\newtheorem*{theoremaux}{Theorem \theoremauxnum}
\gdef\theoremauxnum{1}

\newtheorem{lemma}{\bf Lemma}[section]
 
\newtheorem{theorem}{\bf Theorem}[section]

\newtheorem{corollary}[lemma]{\bf Corollary}
\newtheorem{definition}{\bf Definition}[section]

\journal{~}

\begin{document}

\begin{frontmatter}



\title{A Family of Tetravalent Half-transitive Graphs}



\author[1]{Sucharita Biswas}
\ead{biswas.sucharita56@gmail.com}
\author[2]{Angsuman Das\corref{cor1}}
\ead{angsuman.maths@presiuniv.ac.in}

\address{Department of Mathematics, Presidency University, Kolkata, India}
\cortext[cor1]{Corresponding author}

\begin{abstract}
	In this paper, we introduce a new family of graphs, $\Gamma(n,a)$. We show that it is an infinite family of tetravalent half-transitive Cayley graphs. Apart from that, we determine some structural properties of $\Gamma(n,a)$. 
\end{abstract}

\end{frontmatter}

\section{Introduction}
A graph $G=(V,E)$ is said to be vertex-transitive, edge-transitive and arc-transitive if the automorphism group of $G$, $\mathsf{Aut}(G)$, acts transitively on the vertices, on the edges and on the arcs of $G$ respectively. It is known that an arc-transitive graph is both vertex-transitive and edge-transitive. However, a graph which is both vertex-transitive and edge-transitive may not be arc-transitive, the smallest example being the Holt graph \cite{holt} on $27$ vertices. Such graphs are called {\it half-transitive} graphs. For other definitions related to algebraic graph theory, one is referred to \cite{godsil-royle}.

The study of half-transitive graphs was initiated by Tutte \cite{tutte}, who proved that any half-transitive graph is of even degree. Since any connected $2$-regular is a cycle and a cycle is arc-transitive, the first possibility of finding a half-transitive graph is a $4$-regular or tetravalent graph. The first example of tetravalent half-transitive graphs were given by Bouwer \cite{bouwer} and the smallest example was given in Holt \cite{holt}. Though numerous papers have been published in the last 50 years, the classification of tetravalent half-transitive graphs is not yet complete. In absence of complete classification, two major approaches have been fruitful so far: the first is to characterize half-transitive graphs of some particular orders like $p^3,p^4,p^5,pq,2pq$ etc \cite{char-4},\cite{char-3},\cite{char-2},\cite{char-1}, and the second is to come up with infinite families of half-transitive graphs \cite{chen-li-seress-family},\cite{family-1}.

In this paper, we introduce a new family of tetravalent half-transitive graphs, $\Gamma(n,a)$. It turns out that it is family of Cayley graphs of order $3n$ and for $(n,a)=(9,4)$, we get the Holt graph. In fact, our construction is a generalization of an alternative construction of Holt graph (See the last paragraph of \cite{holt}).

{\definition Let $n$ be a positive integer such that $3\mid \varphi(n)$, where $\varphi$ denote the Euler Totient function. Then $\mathbb{Z}^*_n$, the group of units of $\mathbb{Z}_n$,  is a group of order a multiple of $3$. Let $a$ be a element of order $3$ in $\mathbb{Z}^*_n$ and $b=a^2~(mod ~n)$. Define $\Gamma(n,a)$ to be the graph with vertex-set $\mathbb{Z}_n\times \mathbb{Z}_3$ and the edge-set comprises of edges of the form $(i,j)\sim (ai\pm 1,j-1)$ and $(i,j)\sim (bi\pm b,j+1)$, where the operations in the first and second coordinates are done in modulo $n$ and modulo $3$ respectively.  }

It is obvious that $\Gamma(n,a)$ is tetravalent. One can check that $\Gamma(9,4)$ is the Holt graph. It is also to be noted that for a particular $n$, we can have two graphs, $\Gamma(n,a)$ and $\Gamma(n,b)$. However, these two graphs are isomorphic via the automorphism $\tau: \Gamma(n,a)\rightarrow \Gamma(n,a^2)$ defined by $\tau(i,j)=(ai,-j)$. So, without loss of generality, we assume that $a<b$, where $a,b \in \{2,\ldots,n-2\}$.

On the other hand, let $n$ be a positive integer such that $a_1,b_1,a_2,b_2$ be four elements of order $3$ in $\mathbb{Z}^*_n$ with $a_1b_1\equiv 1~(mod~n)$ and $a_2b_2\equiv 1~(mod~n)$. Then, by the above argument, $\Gamma(n,a_1)\cong\Gamma(n,b_1)$ and $\Gamma(n,a_2)\cong\Gamma(n,b_2)$. However, $\Gamma(n,a_1)$ may not be isomorphic to $\Gamma(n,a_2)$, e.g., for $n=63$, we have $4\cdot 16\equiv 1~(mod ~63)$ and $22\cdot 43\equiv 1~(mod ~63)$, but $\Gamma(63,4)$ is not isomorphic to $\Gamma(63,22)$, as odd girth of $\Gamma(63,4)$ is $9$, whereas that of $\Gamma(63,22)$ is $21$.

The definition of $\Gamma(n,a)$ requires that $3|\varphi(n)$. We discuss the form of $n$, for which this holds. Let $n={p_1}^{\alpha_1}{p_2}^{\alpha_2}\cdots {p_k}^{\alpha_k}$, wher $p_i$ are primes. Then $\varphi(n)={p_1}^{\alpha_1-1}{p_2}^{\alpha_2-1}\cdots {p_k}^{\alpha_k-1}(p_1-1)(p_2-1)\cdots (p_k-1)$. As $3|\varphi(n)$, either $3|{p_i}^{\alpha_i-1}$ or $3|p_i-1$ for some $i$, i.e., $9|n$ or $p_i\equiv 1~(mod~3)$ for some $i$. Thus $n$ is either of the form $9t$ or $pt$ where $p$ is a prime of the form $1~(mod~3)$ and $t$ is a positive integer.  

At this junction, it is important to mention that authors in \cite{alspach-marusic} constructed a family of graphs $M(a;m,n)$ and proved the following theorem.

{\theorem[\cite{alspach-marusic}, Theorem 3.3] \label{Marusic-theorem} Let $n\geq 9$ be odd and $a^3\equiv 1 ~(mod~n)$. Then $M(a;3,n)$ is half-transitive.} 

While the current work was in progress, it caught our attention that our construction $\Gamma(n,a)$, coincidentally, is same as $M(a;3,n)$ (Although the definitions are given in completely different way).  However for two reasons we decided to carry on with the project:
\begin{enumerate}
	\item Firstly, the proof techniques are entirely different: While their proof is built on semiregular automorphisms and blocks, ours is based on $6$-cycles present in the graph.
	\item Secondly and most importantly, we prove that $\Gamma(n,a)$ is half-transitive for all $n$ except $7$ and $14$, i.e., $n$ is not necessarily odd (See Theorem \ref{Marusic-theorem}). In other words, comparing with Theorem \ref{Marusic-theorem}  we prove the half-transitivity of a larger family of graphs.
\end{enumerate} 

In Section \ref{Structure-Section}, we discuss some structural properties of $\Gamma(n,a)$ and in Section \ref{Automorphism-Section}, we discuss the parameters related to the automorphism group of $\Gamma(n,a)$. In Section Appendix, we provide the SageMath \cite{sagemath} code for computing the automorphism group of $\Gamma(n,a)$.

\section{Structural Properties of $\Gamma(n,a)$}\label{Structure-Section}

{\theorem \label{bipartite-theorem} $\Gamma(n,a)$ is bipartite if and only if $n$ is even. }\\
\pf If $n$ is even, then consider the sets $X=\{(i,j): i \mbox{ is even }\}$
and $Y=\{(i,j): i \mbox{ is odd}\}$. Note that as $n$ is even and $a,b \in \mathbb{Z}^*_n$, therefore $a,b$ are odd. Thus, if $i$ is even, then $ai\pm 1$ and $bi\pm b$ are odd, and if $i$ is odd, then $ai\pm 1$ and $bi\pm b$ are even. Thus, no two vertices in $X$ (or $Y$) are adjacent and hence $\Gamma(n,a)$ is bipartite.

On the other hand, if $n$ is odd, as $\Gamma(n,a)$ is Hamiltonian (See Corollary \ref{n=odd-hamiltonian}), it has an odd cycle of length $3n$. Hence $\Gamma(n,a)$ is not bipartite.\qed

{\theorem \label{chromatic-theorem} Then chromatic number $\chi$ of $\Gamma(n,a)$ is given by $$\chi=\left\lbrace \begin{array}{ll}
	2, & \mbox{if }n \mbox{ is even}\\
	3, & \mbox{if }n \mbox{ is odd}
	\end{array}\right.$$  }
\pf The proof is obvious, if $n$ is even. If $n$ is odd, consider the sets $A_0=\{(i,0):i \in \mathbb{Z}_n\}$, $A_1=\{(i,1):i \in \mathbb{Z}_n\}$ and $A_2=\{(i,2):i \in \mathbb{Z}_n\}$. The theorem is true as $A_0,A_1$ and $A_2$ are independent sets in $\Gamma(n,a)$ and $\Gamma(n,a)$ is not bipartite.\qed

{\lemma If $n=p$, where $p\equiv 1~ (mod~3)$ is a prime, then girth of $\Gamma(n,a)$ is $3$.}\\
\pf As $a^3\equiv 1~(mod~n)$ and $a\not\equiv 1~(mod~n)$, we have $n|a^3-1=(a-1)(a^2+a+1)$. As $n$ is prime, we have $a^2+a+1\equiv 0~(mod~n)$. The lemma follows from the observation that $(0,0)\sim (1,2)\sim (a+1,1)\sim (a^2+a+1,0)=(0,0)$ is a $3$-cycle.\qed

{\lemma \label{existence-of-6-cycle} $\Gamma(n,a)$ always have a cycle of length $6$ and hence girth of $\Gamma(n,a)$ is less than or equal to $6$.}\\
\pf  Consider the cycle $(0,0)\sim (1,2)\sim (a-1,1)\sim (a^2-a+1,0)\sim (-a^2+a,2)\sim (b,1)\sim (0,0)$. Clearly it is a $6$-cycle, provided the vertices are distinct. If two vertices are not distinct, then we have either $a-1\equiv b$ or $-a^2+a\equiv 1$ or $a^2-a+1\equiv 0$, i.e., in all cases, we should have $a^2-a+1\equiv 0~(mod~n)$, i.e., $a^3+1\equiv 0~(mod~n)$. This implies $2 \equiv 0~(mod~n)$, a contradiction.\qed

{\lemma \label{no-4-cycle} $\Gamma(n,a)$ does not have any $4$-cycle.}\\
\pf Let, if possible, $\Gamma(n,a)$ has a $4$-cycle. As $\Gamma(n,a)$ is edge-transitive (by Theorem \ref{edge-transitive-theorem}), without loss of generality, we can assume that it has $(0,0)$ and $(1,2)$ as two of its adjacent vertices. Now, the other three neighbours of $(0,0)$ are $X=\{(b,1),(-b,1),(-1,2)\}$ and that of $(1,2)$ are $Y=\{(2b,0),(a+1,1),(a-1,1)\}$. To reach a contradiction, it suffices to show that $X$ and $Y$ does not have any edge between them. 

We start with $(2b,0)\in Y$. If $(2b,0)\sim (-1,2)$, then by the adjacency criterion, we have $2ab\pm 1\equiv -1~(mod~n)$, i.e., $-1\equiv 3,1~(mod~n)$, i.e., $n$ divides $2$ or $4$, a contradiction. If $(2b,0)\sim (b,1)$, then by the adjacency criterion, we have $ab\pm 1\equiv2b~(mod~n)$, i.e., $2b\equiv 0,2~(mod~n)$. Now, $2b\equiv0~(mod~n)$ implies that $b$ is zero divisor in $\mathbb{Z}_n$, a contradiction. So, let $2b\equiv 2~(mod~n)$. If $n$ is odd, this implies $b\equiv 1~(mod~n)$, a contradiction. If $n=2m$ is even, then this implies $b\equiv 1(~mod~m)$, i.e., $b=m+1$ (as $b\neq 1$). However, this in turn implies $b^3=m(m^2+m+3)+1$. Note that $m^2+m+3$ is odd irrespective of $m$ is odd/even. Thus $b^3\neq 1~(mod~2m)$, a contradiction. Thus, $(2b,0)\not\sim (b,1)$. Similarly, it can be shown that $(2b,0)\not\sim (-b,1)$. Thus $(2b,0)$ is not adjacent to any vertices in $X$.

Now, let us consider $(a+1,1)\in Y$. The only possible neighbour of $(a+1,1)$ in $X$ is $(-1,2)$. If $(a+1,1)\sim(-1,2)$, then we have $(a+1)b\pm b\equiv -1~(mod~n)$. As $-1\not\equiv 1~(mod~n)$, we have $1+2b\equiv -1~(mod~n)$, i.e., $2b\equiv -2~(mod~n)$. But, as in previous case, this leads to a contradiction, thereby establishing that $(a+1,1)\not\sim(-1,2)$.

Similarly, it can be shown that $(a-1,1)\in Y$ has no neighbour in $X$. Hence the lemma follows.\qed 

{\corollary \label{n=even-girth=6} If $n$ is even, then girth of $\Gamma(n,a)$ is $6$.}\\
\pf This follows from Theorem \ref{bipartite-theorem}, Lemma \ref{existence-of-6-cycle} and Lemma \ref{no-4-cycle}.

{\lemma \label{9|n-triangle-free} If $9|n$, then $\Gamma(n,a)$ is triangle-free.}\\
\pf If possible, let $\Gamma(n,a)$ conatains a triangle. As $\Gamma(n,a)$ is edge-transitive, without loss of generality, we can assume that two of the vertices of the triangle are $(0,0)$ and $(1,2)$. Then the third vertex is of the form $(x,1)$. As it is adjacent to $(0,0)$, $x\equiv \pm b~(mod~n)$ and as it is adjacent to $(1,2)$, $x\equiv a\pm 1~(mod~n)$. Combining these two, we get $x\equiv a\pm 1\equiv \pm b\equiv a^2~(mod~n)$. Thus, we have four cases:

{\bf Case 1:} $a-1\equiv a^2~(mod~n)$. In this case, multiplying both sides by $a$, we have $a^2-a\equiv1~(mod~n)$, i.e., $n|2$, a contradiction.

{\bf Case 2:} $a+1\equiv a^2~(mod~n)$. In this case, multiplying both sides by $a$, we have
$a^2+a\equiv 1~(mod~n)$, i.e., $2a\equiv 0~(mod~n)$, a contradiction, as $a$ is not a zero-divisor.

{\bf Case 3:} $a-1\equiv -a^2~(mod~n)$. In this case, multiplying both sides by $a$, we have $a^2-a\equiv -1~(mod~n)$, i.e., $2a\equiv 2~(mod~n)$. Again, as in the proof of Lemma \ref{no-4-cycle}, we reach a contradiction.

{\bf Case 4:} $a+1\equiv -a^2~(mod~n)$, i.e., $a^2+a+1 \equiv0~(mod~n)$. As $9|n$, we have $9|(a^2+a+1)$. Moreover, as $a$ is a unit in $\mathbb{Z}_n$, we have $gcd(9,a)=1$. Thus $a$ is of the form $9k+i$, where $i=1,2,4,5,7,8$. Now, it can be checked for all such choices of $i$, $a^2+a+1$ is not a multiple of $9$, a contradiction.

Combining all the cases, the lemma holds.\qed

{\lemma Let $9|n$. Then girth of $\Gamma(n,a)=\left\lbrace \begin{array}{ll}
	5 & \mbox{if }n=9\\
	6 & \mbox{if }n\neq 9
	\end{array} \right.$}\\
\pf If $n$ is an even multiple of $9$, then the lemma follows from  Corollary \ref{n=even-girth=6}. So, we assume that $n=9k$ where $k$ is odd. Clearly $(0,0)\sim (1,2)\sim (5,1)\sim (6,2)\sim (7,1)\sim (0,0)$ is a $5$-cycle in $\Gamma(9,4)$, the Holt graph. Thus from Lemma \ref{no-4-cycle} and Lemma \ref{9|n-triangle-free}, it follows that girth of $\Gamma(9,4)$ is $5$. 

If possible, let $n=9k$ where $k\geq 3$ is an odd number and girth of $\Gamma(n,a)$ is $5$. Let $C$ be a $5$-cycle in $\Gamma(n,a)$. As $\Gamma(n,a)$ is edge-transitive (by Theorem \ref{edge-transitive-theorem}), without loss of generality, we can assume that $C$ has $(0,0)$ and $(1,2)$ as two of its adjacent vertices. Let $X=\{(-1,2),(b,1),(-b,1)\}$ be the other three neighbours of $(0,0)$ and $Y=\{(2a^2,0),(a+1,1),(a-1,1)\}$ be the other three neighbours of $(1,2)$. Since $C$ is a $5$-cycle, some vertex in $X$ and some vertex in $Y$ must have a common neighbour.

We will show that this is not possible. As each $X$ and $Y$ has $3$ vertices, we need to check that all possible $9$ pairs  do not have any common neighbour. However, as the proof technique is similar for all the pairs, we show it only for one pair, namely $(b,1)=(a^2,1)$ and $(a+1,1)$. The common neighbour, if any, of $(a^2,1)$ and $(a+1,1)$, is of the form $(x,0)$ or $(x,2)$.

{\bf Case 1:} Let $(x,0)$ be the common neighbour of $(a^2,1)$ and $(a+1,1)$. Then, by the adjacency criterion, we have $x\equiv a^3\pm 1\equiv a^2+a\pm 1~(mod~n)$, i.e., $a^2+a\pm 1\equiv 0$ or $2~(mod~n)$. If $a^2+a\pm 1\equiv 0~(mod`n)$, then we land in Case 3 or Case 4 of the proof of Lemma \ref*{9|n-triangle-free} and hence a contradiction. If $a^2+a+ 1\equiv 2~(mod`n)$, then also we land in Case 3 of Lemma \ref*{9|n-triangle-free} and hence a contradiction. If $a^2+a- 1\equiv 2~(mod`n)$, i.e., $a^2+a- 3\equiv 0~(mod`n)$, then multiplying both sides by $a$, we get $a^2-3a+1\equiv 0~(mod~n)$. Hence, subtracting on eequation from the other, we have $2(1-a)\equiv 0~(mod ~n)$. Now, as $n$ is odd, this implies $a\equiv 1~(mod~n)$, a contradiction to the fact that $\circ(a)=3$ in $\mathbb{Z}_n^*$. 

{\bf Case 2:} Let $(x,1)$ be the common neighbour of $(a^2,1)$ and $(a+1,1)$. Then, by the adjacency criterion, we have $x\equiv 1+a^2\pm a^2\equiv a\pm a^2~(mod~n)$, i.e., $a\pm a^2\equiv 1$ or $1+2a^2~(mod~n)$. If $a\pm a^2\equiv 1$, then we reach either Case 1 or Case 3 of the proof of Lemma \ref*{9|n-triangle-free} and hence a contradiction. Also, if $a+a^2\equiv 1+2a^2~(mod~n)$, we reach Case 1 of the proof of Lemma \ref*{9|n-triangle-free} and hence a contradiction. If $a-a^2\equiv 1+2a^2~(mod~n)$, i.e., $3a^2-a+1\equiv 0~(mod~n)$, then multiplying both sides by $a$, we get $a^2-a-3\equiv 0~(mod~n)$. Subtracting one equation from the other yields
$2(2a+1)\equiv 0~(mod~n)$. As $n$ is odd, this implies $2a\equiv -1~(mod~n)$. Cubing both the sides, we get $8\equiv 8a^3\equiv -1~(mod~n)$, which is possible only if $n=9$, a contradiction.

Finally, combining both the cases, the lemma holds.\qed

\section{Automorphisms of $\Gamma(n,a)$}\label{Automorphism-Section}
Let $G=\mathsf{Aut}(\Gamma(n,a))$. We start by noting the following automorphisms of $\Gamma(n,a)$. 
$$\begin{array}{lll}
\alpha: (i,j) \mapsto (i+a^{-j},j); & ~~\beta : (i,j)\mapsto (i,j+1); & ~~\gamma: (i,j) \mapsto (-i,j);
\end{array}$$
It can be shown that $\alpha,\beta,\gamma \in G$ and $\circ(\alpha)=n,\circ(\beta)=3$ and $\circ(\gamma)=2$. Moreover, we have the following relations: $\alpha\beta=\beta\alpha^{a^2}, \alpha\gamma=\gamma\alpha^{-1}$ and $\beta\gamma=\gamma\beta$.

{\theorem \label{cayley-theorem} $\Gamma(n,a)$ is a Cayley graph.}\\
\pf Let $H=\langle \alpha, \beta \rangle$. Clearly it forms a subgroup of $G$. Also as $\circ(\alpha)=n,\circ(\beta)=3$ and $\alpha\beta=\beta\alpha^{a^2}$, we have $$H=\{\alpha^i\beta^j:0\leq i \leq n-1,0\leq j \leq 2\} \mbox{ and }|H|=3n=|\Gamma(n,a)|.$$
We will show that $H$ acts regularly on $\Gamma(n,a)$. As $|H|=|\Gamma(n,a)|$, it is enough to show that $H$ acts transitively on $\Gamma(n,a)$. As $i\mapsto i+a^{-j}$ is a permutation of $\mathbb{Z}_n$ order $n$ and $j\mapsto j+1$ is a permutation of $\mathbb{Z}_3$ order $3$, the action of $H$ on $\Gamma(n,a)$ is transitive. \qed

Note that $H$ is a semidirect product of $\langle \alpha \rangle$ and $\langle \beta \rangle$, as $\beta^{-1}\alpha\beta=\alpha^{a^2}$ and $a^2$ and $n$ are coprime, and $\Gamma(n,a)=Cay(H,S)$ where $S=\{\beta^2\alpha,\beta^2\alpha^{-1},\beta\alpha^b,\beta\alpha^{-b} \}$. We now recall a result on hamiltonicity of Cayley graphs.

{\theorem[\cite{marusic-hamiltonian}, Theorem 3.3] \label{marusic-ham-theorem} Every connected Cayley graph of a semidirect product of a cyclic group of prime order by an abelain group of odd order is Hamiltonian. \qed}

{\corollary \label{n=odd-hamiltonian} $\Gamma(n,a)$ is Hamiltonian, if $n$ is odd.}\\
\pf By Theorem \ref{cayley-theorem}, we have $\Gamma(n,a)$ is a Cayley graph on $H$ and $H$ is a semidirect product of cyclic group of order $3$, namely $\langle \beta \rangle$ and another cyclic group of odd order $n$, namely $\langle \alpha \rangle$. Thus, by Theorem \ref{marusic-ham-theorem}, $\Gamma(n,a)$ is Hamiltonian.\qed

{\theorem \label{edge-transitive-theorem} $\Gamma(n,a)$ is edge-transitive.}\\
\pf As $\Gamma(n,a)$ is Cayley, it is vertex-transitive. Hence, it is enough to show that any two edges incident to $(0,0)$ can be permuted by an automorphism. As $\Gamma(n,a)$ is tetravalent, the four vertices adjacent to $(0,0)$ are namely: $(1,2),(-1,2),(b,1)$ and $(-b,1)$. Let us name the following edges: $$\begin{array}{rlcrl}
e_1: & (0,0)\sim (1,2)& ~~ &e_2: & (0,0)\sim (-1,2)\\ 
e_3: & (0,0)\sim (b,1)& ~~ &e_4: & (0,0)\sim (-b,1)
\end{array}$$
It is to be noted that $\gamma(e_1)=e_2$, $\alpha\beta\gamma(e_1)=~\stackrel{\leftarrow}{e_3}$ and $\gamma\alpha\beta\gamma(e_1)=~\stackrel{\leftarrow}{e_4}$. The reverse arrow on top denote that the orientation of the edge changed. Hence, the theorem.\qed 

For $n=7,14$, SageMath \cite{sagemath} computation shows that $\Gamma(n,a)$ is arc-transitive. Next we prove that $\Gamma(n,a)$ is not arc-transitive if $n\neq 7,14$. For that, we show that there does not exist any graph automorphism $\varphi$ which maps the arc $e_3$ to $e_1$, i.e.,  $\varphi((0,0))=(0,0)$ and $\varphi((b,1))=(1,2)$. 

If possible, let such an automorphism $\varphi$ exists with $\varphi((0,0))=(0,0)$ and $\varphi((b,1))=(1,2)$. As $(1,2)$ is adjacent to $(0,0)$, its image under $\varphi$ should be one among $(b,1)$, $(-b,1)$ and $(-1,2)$. However, the next theorem shows that this can not hold.

{\theorem \label{main-theorem} If $\varphi$ is an automorphism of $\Gamma(n,a)$ such that $n\neq 7,14$ and $\varphi((0,0)=(0,0)$ and $\varphi((b,1))=(1,2)$, then $\varphi((1,2)) \notin \{(b,1),(-b,1),(-1,2)\}$.\qed}

Thus, from Theorem \ref{cayley-theorem}, Theorem \ref{edge-transitive-theorem} and Theorem \ref{main-theorem}, we get the desired result as follows:

{\theorem  $\Gamma(n,a)$ is half-transitive if $n\neq 7,14$. \qed}

\section{Proof of Theorem \ref{main-theorem} }
To prove Theorem \ref{main-theorem}, we prove a lemma and three theorems. Throughout this section, $\varphi$ denote an automorphism of $\Gamma(n,a)$ and $G$ denote the full automorphsim group of $\Gamma(n,a)$.

{\lemma \label{relations} The following relations can not  hold.\\
	\textbf{1.} $2a-4b=0$, \textbf{2.} $2a+4b=0$ except $n=9$, \textbf{3.} $4a-2b=0$ except $n=7,14$, \textbf{4.} $4a+2b=0$ except $n=18$,  \textbf{5.} $2a-2b=0$, \textbf{6.} $2a+2b=0$, \textbf{7.} $4a+4=0$, \textbf{8.} $2a+6=0$,  \textbf{9.} $2(a+b-1)=0$, \textbf{10.} $2(b-a+1)=0$, \textbf{11.} $2(a-b+1)=0$, \textbf{12.} $2(a+b+2)=0$,  \textbf{13.} $2(a-b-2)=0$.  }\\
	\pf \textbf{1.} $2a-4b$, i.e., $8=64$, i.e., $56=0$, i.e., $n \mid 56$, hence $n=7, 14, 28, 56$ and the possible values of $b$ are $4,11,25,25$ respectively. In all these cases $2b \neq 4$, which is a contradiction.\\
	\textbf{2.} $2a+4b=0$, i.e., $8=-64$, i.e., $72=0$, i.e., $n \mid 72$, hence $n= 9, 18, 36, 72$ and the possible values of $b$ are $7,13,25,49$. But the relation holds only if $n=9$ and $b=7$. \\
	\textbf{3.}  $4a-2b=0$, $8=64$, i.e., $56=0$, i.e., $n \mid 56$, hence $n=7, 14, 28, 56$ and the possible values of $(a,b)$ are $(2,4),(9,11),(9,25),(9,25)$ respectively. But the relation holds only if $n=7,14$.\\
	\textbf{4.} $4a+2b=0$, i.e., $64=-8$, i.e., $72=0$, i.e., $n \mid 72$, hence $n= 9, 18, 36, 72$ and the possible values of $(a,b)$ are $(4,7),(7,13),(13,25),(25,49)$ respectively. But the relation holds only if $n=18$.\\
	\textbf{5.} $2a-2b=0$, i.e., $2a\equiv 2(mod~n)$. If $n$ is odd then $a=1$, which is impossible. Let $n$ be even and $n=2m$. Then we have $m \mid a-1$, i.e., $a=mt+1$, for some $t \in \Bbb Z$. As  $a \neq 1$ so $a=m+1$, i.e., $a^3-1= m(m^2+3m+3)$. Note that irrespective of $m$ is odd or even, $(m^2+3m+3)$ is odd, say $(2s+1)$, for some $s \in \Bbb Z$. So we have $a^3-1=m(2s+1)$, i.e., $a^3-1 \equiv m (mod~n)$, which is a contradiction.  \\
	\textbf{6.} The proof is same as \textbf{5.}   \\
	\textbf{7.} $4a+4=0$, i.e., $4a\equiv -4(mod~n)$. If $n$ is odd then $a=-1$, which is impossible. Let $n$ is even and $n=2m$, then we have $m \mid 2(a+1)$, i.e., $2a=mt-2$, for some $t \in \Bbb Z$. As $2a \neq -2$ so  $2a=m-2$, i.e., $8(a^3-1)=m(m^2-6m+12)-16$. If $m$ is even then $(m^2-6m+12)$ is even, say $2u$, for some $u \in \Bbb Z$. So we have  $8(a^3-1)=2mu-16$, i.e., $8(a^3-1) \equiv -16(mod~n)$, which is a contradiction. If $m$ is odd then $(m^2-6m+12)$ is odd, say $2v+1$, for some $v \in \Bbb Z$.  So we have  $8(a^3-1)=m(2v+1)-16$, i.e., $8(a^3-1) \equiv m-16(mod~n)$, which is a contradiction as $m \neq 16$.\\
	\textbf{8.} $2a+6=0$, i.e., $8=-216$, i.e., $224=0$, i.e., $n \mid 224$, i.e., $n= 7, 14, 28, 56, 112, 224$. However, in all these cases, the possible values of $a$ does not allow $2a+6=0$.\\
	\textbf{9.} $2(a+b-1)=0$, i.e., $2(1+a-b)=0$, i.e., $4a=0$, contradicting that $a$ is a unit.\\
	\textbf{10.} The proof is same as \textbf{9}.\\
	\textbf{11.} $2(a-b+1)=0$, i.e., $2(1-a+b)=0$, i.e., $2(a-b+1)+2(1-a+b)=0$, i.e.,  $4=0$, which is a contradiction.\\
	\textbf{12.} $2(a+b+2)=0$, i.e., $2(1+a+2b)=0$, i.e., $4(a+b+2)-2(1+a+2b)=0$, i.e., $2(a+3)=0$, i.e.,  $8=-216$, i.e., $224=0$, i.e., $n \mid 224$, i.e., $n= 7, 14, 28, 56, 112, 224$. In all these cases $2a+6 \neq 0$, which is a contradiction. \\
	\textbf{13.} $2(a-b-2)=0$, i.e., $2(1-a-2b)=0$, i.e., $2(a-b-2)+2(1-a-2b)=0$, i.e., $6b+2=0$, i.e., $2a+6=0$. Rest of the proof is same as \textbf{$8$}.\qed

{\theorem \label{phi(1,2)=(-1,2)}  If $\varphi \in G$ and $\varphi((0,0))=(0,0),~ \varphi((b,1))=(1,2)$, $\varphi((1,2))=(-1,2)$ then $n=7$ or $14$. }\\
\pf Consider the cycle $C: (0,0)\sim (b,1) \sim (a+b,2) \sim (1+a+b,0)\sim (a+1,1)\sim (1,2)\sim (0,0)$. Then $\varphi(C):(0,0)\sim (1,2) \sim \varphi((a+b,2)) \sim \varphi((1+a+b,0)) \sim \varphi((a+1,1))\sim (-1,2)\sim (0,0)$. As $\varphi((a+b,2)) \sim (1,2)$ and $\varphi((0,0))=(0,0)$ so $\varphi((a+b,2)) \in  \lbrace (2b,0), (a \pm 1,1) \rbrace$.  Again, since $\varphi((a+1,1)) \sim (-1,2)$ and $\varphi((0,0))=(0,0)$ imply $\varphi((a+1,1)) \in  \lbrace (-2b,0), (-a \pm 1,1) \rbrace$. Also $\varphi((1+a+b,0)) \sim \varphi((a+b,2))$ and $\varphi((b,1))=(1,2)$ imply 
\begin{equation}\label{image-of-(1+a+b,0)00}
\varphi((1+a+b,0)) \in \lbrace (2a \pm b,1), (3,2), (1+ 2b,2), (b +a \pm 1, 0),(1-2b,2), (b-a \pm 1,0)  \rbrace.
\end{equation}

If $\varphi((a+1,1))=(-2b,0)$ then $\varphi((1+a+b,0)) \sim \varphi((a+1,1))$ and $\varphi((1,2))=(-1,2))$ imply
\begin{equation}\label{image-of-(1+a+b,0)01}
\varphi((1+a+b,0)) \in \lbrace (-3,2), (-2a \pm b,1) \rbrace 
\end{equation}
From the Equations \ref{image-of-(1+a+b,0)00} and \ref{image-of-(1+a+b,0)01} we have, 
\begin{itemize}
\item either $-3=3$, i.e., $6=0$, i.e., $n=6$ which is impossible. 
\item or $-3=1+2b$, i.e., $2a+4b=0$, which is possible only when $n=9$, (by Lemma \ref{relations}). However, direct Sagemath computation for $n=9$ shows that such $\varphi$ does not exist.
\item or $-3=1-2b$, i.e., $2a-4b=0$,  which is impossible by Lemma \ref{relations}.
\item or $-2a \pm b=2a \pm b$, i.e., $4a=0$ or $4a-2b=0$ or $4a+2b=0$. Though the first one is impossible, the other two can hold only if $n=7,14,18$ (by Lemma \ref{relations}). However, direct Sagemath computation for $n=7,14$ and $18$ shows that such $\varphi$ does not exist. 	
\end{itemize}
Hence $\varphi((a+1,1)) \neq (-2b,0)$.

If $\varphi((a+1,1))=(-a+1,1)$, then $\varphi((1+a+b,0)) \sim \varphi((a+1,1))$ and $\varphi((1,2))=(-1,2)$ imply
\begin{equation}\label{image-of-(1+a+b,0)02}
\varphi((1+a+b,0)) \in \lbrace (-1+2b,2), (-b+a \pm 1,0) \rbrace 
\end{equation}
From the Equations \ref{image-of-(1+a+b,0)00} and \ref{image-of-(1+a+b,0)02} we have
\begin{itemize}
\item either $-b+a \pm 1= b-a \pm 1$, i.e., $2(a-b)=0$ or $2(b-a+1)=0$ or $2(b-a-1)=0$, all of which are impossible by Lemma \ref{relations}. 
\item or $-b+a \pm 1=b+a \pm 1$, i.e., $2b=0$ or$2a-2b=0$ or $2a+2b=0$,  all of which are impossible by Lemma \ref{relations}.
\item or $-1+2b=3$, i.e., $2a-4b=0$,  which is impossible by Lemma \ref{relations}.
\item or $-1+2b=1+2b$, i.e., $2=0$, which is a contradiction.
\item or $-1+2b=1-2b$, i.e., $4a-2b=0$ which can hold only if $n=7$ or $14$. (by Lemma \ref{relations}.) However, direct Sagemath computation for $n=7,14$ shows that such $\varphi$ does not exist.
\end{itemize}  
 Hence $\varphi((a+1,1)) \neq (-a+1,1)$.

If $\varphi((a+1,1))=(-a-1,1)$ then $\varphi((1+a+b,0)) \sim \varphi((a+1,1))$ and $\varphi((1,2))=(-1,2))$ imply
\begin{equation}\label{image-of-(1+a+b,0)03}
\varphi((1+a+b,0)) \in \lbrace (-1-2b,2), (-b-a \pm 1,0) \rbrace 
\end{equation}
From the Equations \ref{image-of-(1+a+b,0)00} and \ref{image-of-(1+a+b,0)03} we have 
\begin{itemize}
\item either $-1-2b=3$, i.e., $2a+4b=0$,  which can hold only if $n=9$. (by Lemma \ref{relations}). However, direct Sagemath computation rules out this possibility.
\item or $-1-2b=1+2b$, i.e., $4a+2b=0$, which can hold only if $n=18$. (by Lemma \ref{relations}). However, direct Sagemath computation rules out this possibility.
\item or $-1-2b=1-2b$, i.e., $2=0$,  which is a contradiction.
\item or $-b-a \pm 1= b-a \pm 1$, i.e., $2b=0$, or $2a+2b=0$, or $2a-2b=0$  all of which are impossible by Lemma \ref{relations}.
\item or $-b-a \pm 1=b+a \pm 1$, i.e., $2(b+a-1)=0$ or $2(b+a)=0$  (which are impossible by Lemma \ref{relations}), but $2(1+a+b)=0$ may hold.
\end{itemize} 

Therefore we have $\varphi((1+a+b,0))=(1+a+b,0)$, $\varphi((a+1,1))=(-a-1,1)$, $\varphi((a+b,2))=(a+1,1)$ with $2(a+b+1)=0$. 

Consider the cycle $C': (1+a+b,0)\sim (a+1,1) \sim (1+2b,2) \sim (2a,0) \sim (b+2,1) \sim (a+b,2) \sim (1+a+b,0)$. Then $\varphi(C'): (1+a+b,0)\sim (-a-1,1) \sim \varphi((1+2b,2)) \sim \varphi((2a,0)) \sim \varphi((b+2,1)) \sim (a+1,1) \sim (1+a+b,0)$. Now $\varphi((b+2,1)) \sim (a+1,1)$,  $\varphi((1+a+b,0))=(1+a+b,0)$ and  $\varphi((b,1))=(1,2)$ imply $\varphi((b+2,1)) \in \lbrace (1+2b,2), (b+a-1,0)\rbrace$.  Again $\varphi((1+2b,2)) \sim (-a-1,1)$, $\varphi((a+b+1,0))=(a+b+1,0)=(-a-b-1,0)$  and  $\varphi((1,2))=(-1,2)$ imply $\varphi((1+2b,2)) \in \lbrace (-1-2b,2), (-b-a+ 1,0) \rbrace$. Also $\varphi((2a,0)) \sim \varphi((b+2,1))$ and $\varphi((a+b,2))=(a+1,1)$ imply
\begin{equation}\label{image-of-(2a,0)0} 
\varphi((2a,0)) \in \lbrace  (b+2a \pm b,0), (a+3,1), (a+1-2b,1), (1+b-a \pm 1,2) \rbrace.
\end{equation}

Let $\varphi((1+2b,2))=(-1-2b,2)$. $\varphi((2a,0))\sim \varphi((1+2b,2))$ and  $\varphi((a+1,1))=(-a-1,1)$ imply
\begin{equation}\label{image-of-(2a,0)1}
\varphi((2a,0)) \in \lbrace (-b-2a \pm b,0), (-a-3,1) .
\end{equation}
From the Equations \ref{image-of-(2a,0)0} and \ref{image-of-(2a,0)1} we have
\begin{itemize}
	\item either  $-b-2a \pm b=b+2a \pm b$, i.e., $4a+4=0$ or $4a=0$ (which are impossible by Lemma \ref{relations}) or $4a+2b=0$, which can hold only if $n=18$. However, direct Sagemath computation rules out this possibility.
	\item or $-a-3=a+3$, i.e., $2a+6=0$, which is impossible by Lemma \ref{relations}.
	\item or $-a-3=a+1-2b$, i.e., $2(a-b+2)=0$. Also, we had $2(a+b+1)=0$ previously. This yields $2a=4$, i.e., $n=7$ or $14$.
\end{itemize}
 Hence $\varphi((1+2b,2))=(-1-2b,2)$ is possible only if $n=7$ or $14$. Moreover,  direct Sagemath computation for $n=7$ and $14$ confirms the possibility.

Let $\varphi((1+2b,2))=(-b-a+1,0)$. $\varphi((2a,0))\sim \varphi((1+2b,2))$ and  $\varphi((a+1,1))=(-a-1,1)$ imply
\begin{equation}\label{image-of-(2a,0)2}
\varphi((2a,0)) \in \lbrace (-a-1+2b,1), (-1-b+a \pm 1,2) .
\end{equation}
From the Equations \ref{image-of-(2a,0)0} and \ref{image-of-(2a,0)2} we have, 
\begin{itemize}
\item either $-1-b+a \pm 1=1+b-a \pm 1$, i.e., $2(b-a+1)=0$ or $2(a-b)=0$ or $2(a-b-2)=0$, all of which are impossible by Lemma \ref{relations}.
\item or $-a-1+2b=a+1-2b$, i.e., $2(a+1-2b)=0$, i.e., $2(a+b-2)=0$. Hence combining $2(a+b+1)=0$ and $2(a+b-2)=0$, we have $6=0$, which is impossible.
\item or $-a-1+2b=a+3$, i.e., $2(a-b+2)=0$.  Therefore from $2(a+b+1)=0$ and $2(a-b+2)=0$ we have $2a=4$, i.e., $n=7$ or $14$. 
\end{itemize}
Hence $\varphi((1+2b,2))=(-b-a+1,0)$ may be possible if $n=7$ or $14$. Moreover,  direct Sagemath computation for $n=7$ and $14$ confirms the possibility. 

Therefore for  $\varphi \in G$ we can have  $\varphi((0,0))=(0,0),~ \varphi((b,1))=(1,2)$, $\varphi((1,2))=(-1,2)$ only if $n=7$ or $14$.   \qed

{\theorem \label{phi(1,2)=(-b,1)} If $\varphi \in G$ and $\varphi((0,0))=(0,0),~ \varphi((b,1))=(1,2)$ then  $\varphi((1,2))\neq (-b,1)$.  }\\
\pf If possible let  $\varphi \in G$ and $\varphi((0,0))=(0,0),~ \varphi((b,1))=(1,2)$ and $\varphi((1,2))= (-b,1)$.   Consider the cycle $C: (0,0)\sim (b,1) \sim (a+b,2) \sim (1+a+b,0)\sim (a+1,1)\sim (1,2)\sim (0,0)$. Then $\varphi(C):(0,0)\sim (1,2) \sim \varphi((a+b,2)) \sim \varphi((1+a+b,0)) \sim \varphi((a+1,1))\sim (-b,1)\sim (0,0).$ As $\varphi((a+b,2)) \sim (1,2)$ and $\varphi((0,0))=(0,0)$ so $\varphi((a+b,2)) \in  \lbrace (2b,0), (a \pm 1,1) \rbrace$. Again,  $\varphi((a+1,1)) \sim (-b,1)$ and $\varphi((0,0))=(0,0)$ imply $\varphi((a+1,1)) \in  \lbrace (-2,0), (-a \pm b,2) \rbrace$. Now as  $\varphi((1+a+b,0)) \sim \varphi((a+b,2))$ and $\varphi((b,1))=(1,2)$ then we have,
\begin{equation}\label{image-of-(1+a+b,0)10}
\varphi((1+a+b,0)) \in \lbrace (2a \pm b,1), (3,2), (1+ 2b,2), (b +a \pm 1, 0),(1-2b,2), (b-a \pm 1,0)  \rbrace.
\end{equation}

Depending upon the value of $\varphi((a+1,1))$, one of the three cases, namely ({\bf Case A:} $\varphi((a+1,1))=(-2,0)$), ({\bf Case B:} $\varphi((a+1,1))=(-a-b,2)$) and ({\bf Case C:} $\varphi((a+1,1))=(-a+b,2)$) must hold. However, before resolving this three cases, we prove a claim which will be crucial in the following proof.

{\bf Claim: }$\varphi((-1-a-b,0)) \in \lbrace (-2a \pm b,1), (-3,2), (-1+ 2b,2), (-b +a \pm 1, 0),$

\hspace{5.5cm} $(-1-2b,2), (-b-a \pm 1,0),(2a \pm 1,2), (3b,1), (b+2,1)$, 

\hspace{5.5cm} $(1+a \pm b, 0),(b-2,1), (1-a \pm b,0)  \rbrace.$

{\bf Proof of Claim:} As $(-b,1)\sim (0,0)$ and $(-1,2) \sim (0,0)$, we have $\varphi((-b,1)), \varphi((-1,2)) \in \lbrace (b,1), (-1,2) \rbrace$.

\textbf{Case 1:} Let $\varphi((-b,1))=(-1,2)$ and $\varphi((-1,2))=(b,1)$.
Consider the cycle $C': (0,0)\sim (-b,1) \sim (-a-b,2) \sim (-1-a-b,0)\sim (-a-1,1)\sim (-1,2)\sim (0,0)$, then $\varphi(C'):(0,0)\sim (-1,2) \sim \varphi((-a-b,2)) \sim \varphi((-1-a-b,0)) \sim \varphi((-a-1,1))\sim (b,1)\sim (0,0)$. As $\varphi((-a-b,2)) \sim (-1,2)$ and $\varphi((0,0))=(0,0)$ so $\varphi((-a-b,2)) \in  \lbrace (-2b,0), (-a \pm 1,1) \rbrace$.  $\varphi((-a-1,1)) \sim (b,1)$ and $\varphi((0,0))=(0,0)$ imply $\varphi((-a-1,1)) \in  \lbrace (2,0), (a \pm b,2) \rbrace$. Now $\varphi((-1-a-b,0)) \sim \varphi((-a-b,2))$ and $\varphi((-b,1))=(-1,2)$ imply 
\begin{equation}\label{image-of-(-1-a-b,0)0}
\varphi((-1-a-b,0)) \in \lbrace (-2a \pm b,1), (-3,2), (-1+ 2b,2), (-b +a \pm 1, 0),(-1-2b,2), (-b-a \pm 1,0)  \rbrace.
\end{equation}

\textbf{Case 2:} Let $\varphi((-b,1))=(b,1)$ and $\varphi((-1,2))=(-1,2)$.
Consider the cycle $C': (0,0)\sim (-b,1) \sim (-a-b,2) \sim (-1-a-b,0)\sim (-a-1,1)\sim (-1,2)\sim (0,0)$, then $\varphi(C'):(0,0)\sim (b,1) \sim \varphi((-a-b,2)) \sim \varphi((-1-a-b,0)) \sim \varphi((-a-1,1))\sim (b,1)\sim (0,0)$. As $\varphi((-a-b,2)) \sim (b,1)$ and $\varphi((0,0))=(0,0)$ so $\varphi((-a-b,2)) \in  \lbrace (-2,0), (a \pm b,2) \rbrace$.  $\varphi((-a-1,1)) \sim (-1,2)$ and $\varphi((0,0))=(0,0)$ imply $\varphi((-a-1,1)) \in  \lbrace (-2b,0), (-a \pm 1,1) \rbrace$. Now $\varphi((-1-a-b,0)) \sim \varphi((-a-b,2))$ and $\varphi((-b,1))=(b,1)$ imply 
\begin{equation}\label{image-of-(-1-a-b,0)1}
\varphi((-1-a-b,0)) \in \lbrace (2a \pm 1,2), (3b,1), (b+2,1), (1+a \pm b, 0),(b-2,1), (1-a \pm b,0)  \rbrace.
\end{equation}
Combining Case 1 and 2,  the claim follows.

Now, we turn towards the three cases mentioned earlier.

{\bf Case A:} If $\varphi((a+1,1))=(-2,0)$ then $\varphi((1+a+b,0)) \sim \varphi((a+1,1))$ and $\varphi((1,2))=(-b,1)$ imply
\begin{equation}\label{image-of-(1+a+b,0)11}
\varphi((1+a+b,0)) \in \lbrace (-3b,1), (-2a \pm 1,2) \rbrace 
\end{equation}
From the Equations \ref{image-of-(1+a+b,0)10} and \ref{image-of-(1+a+b,0)11} we have, 
\begin{itemize}
	\item either $-3b=2a \pm b$, i.e., $-2b=2a$ or $2b=-4$. By Lemma \ref{relations}, this can hold only if $n=9$. However, direct SageMath computation for $n=9$ show that such $\varphi$ does not exist.
	\item or $-2a \pm 1=3$, i.e., $-2a=4$ or $-2a=2$, i.e., $4a+2b=0$ or $2a+2b=0$. By Lemma \ref{relations}, $2a+2b=0$ can not hold and $4a+2b=0$ can hold only if $n=18$. However, direct SageMath computation for $n=18$ shows that such $\varphi$ does not exist.
	\item or $-2a \pm 1=1-2b$, i.e., $2a=2b$, $2(b-a-1)=0$, both of which are impossible by the Lemma \ref{relations}.
	\item or $-2a \pm 1 = 1 +2b$, i.e., $2a+2b=0$ (which is impossible by the Lemma \ref{relations}) but 
	\begin{equation} \label{2a+2b+2=0,00}
	2a+2b+2=0~ may~ hold.
	\end{equation}
	
\end{itemize}

When $\varphi((a+1,1))=(-2,0)$, $\varphi((1+a+b,0))=(1+2b,2)$, $\varphi((a+b,2))=(a+1,1)$ then we have the Equation \ref{2a+2b+2=0,00}. As $2a+2b+2=0$, i.e., $a+b+1=-a-b-1$, then $\varphi((-1-a-b,0)) = (1+2b,2)$.  But from the Equations \ref{image-of-(-1-a-b,0)0} and \ref{image-of-(-1-a-b,0)1} we have $\varphi((-1-a-b,0)) \neq (1+2b,2)$, which is a contradiction. Hence $\varphi((a+1,1)) \neq (-2,0)$ and {\bf Case A} can not hold.

{\bf Case B:} If $\varphi((a+1,1))=(-a-b,2)$  then $\varphi((1+a+b,0)) \sim \varphi((a+1,1))$ and $\varphi((1,2))=(-b,1))$ imply
\begin{equation}\label{image-of-(1+a+b,0)12}
\varphi((1+a+b,0)) \in \lbrace (-b-2,1), (-1-a\pm b,0) \rbrace 
\end{equation}
From the Equations \ref{image-of-(1+a+b,0)10} and \ref{image-of-(1+a+b,0)12} we have, \\
    \textbf{Case B(1):} either $-b-2=2a \pm b$, i.e., $2a+2b=0$, which is impossible  by the Lemma \ref{relations}, or
	\begin{equation} \label{2a+2b+2=0,01}
	2a+2b+2=0~ may~ hold.
	\end{equation}
	\textbf{Case B(2):} or we have $-1-a\pm b=b+a \pm 1$, i.e., $2a+2b=0$ or $2a=0$ (which are impossible by Lemma \ref{relations}) but  $-1-a-b=b+a+1$, i.e.,
	\begin{equation} \label{2a+2b+2=0,02}
	2a+2b+2=0~ may~ hold.
	\end{equation}
	\textbf{Case B(3):} or we have $-1-a \pm b=b-a \pm 1$. This gives rise to four equations, out of which three are impossible by Lemma \ref{relations}, namely $2=0,~2b=0$ and $2a+2b=0$. The only possibility which remains is $-1-a+b=b-a-1$ and it is an identity.

 So assuming this identity, we have $\varphi((a+1,1))=(-a-b,2)$, $\varphi((a+b+1,0))=(b-a-1,0)$ and $\varphi((a+b,2))=(a-1,1)$. Similarly we can show that $\varphi((a-1,1))=(-a+b,2)$, $\varphi((b-a+1,0))=(b+a-1,0)$ and $\varphi((a-b,2))=(a+1,1)$.  Now $\varphi((a+b,2))=(a-1,1)$, $\varphi(((a-b,2))=(a+1,1)$ and $\varphi((2,0)) \sim \varphi((b,1))=(1,2)$ imply $\varphi((2,0))=(2b,0)$.

Now, consider the cycle $C_2: (a+b+1,0) \sim (a+1,1) \sim (1,2) \sim (2b,0) \sim (2a+b,1) \sim (a+b+2,2) \sim (a+b+1,0).$ So $\varphi(C_2): (b-a-1,0) \sim (-a-b,2) \sim (-b,1) \sim  \varphi((2b,0)) \sim \varphi((2a+b,1)) \sim \varphi((a+b+2,2)) \sim (b-a-1,0).$

 Again $\varphi((0,0))=(0,0)$, $\varphi((a+1,1))=(-a-b,2)$ and $\varphi((2b,0)) \sim (-b,1)$ imply $\varphi((2b,0)) \in \lbrace (-a+b,2), (-2,0) \rbrace$. And $\varphi(((a+b,2))=(a-1,1)$, $\varphi((a+1,1))=(-a-b,2)$ and $\varphi((a+b+2,2)) \sim  (b-a-1,0)$ imply $\varphi((a+b+2,2)) \in \lbrace (-b-a+2,2) , (-1-2b+a,1) \rbrace$. $\varphi((1,2))=(-b,1)$ and $\varphi((2a+b,1)) \sim \varphi((2b,0)) $ imply 
\begin{equation}\label{(a+b+2,2)00}
\varphi((2a+b,1)) \in \lbrace (-1+a \pm b,0), (-b+2,1), (-3b,1), (-2a \pm 1,2).
\end{equation}
\textbf{Case B(3)(a):} If $\varphi((a+b+2,2))=(-b-a+2,2)$ then $\varphi((a+b+1,0))=(b-a-1,0)$ and $\varphi((2a+b,1)) \sim \varphi((a+b+2,2))$ imply 
\begin{equation}\label{(a+b+2,2)01}
\varphi((2a+b,1)) \in \lbrace (a-1+3b,0), (-1-b+2a \pm 1,1) \rbrace.
\end{equation}
From the Equations \ref{(a+b+2,2)00} and \ref{(a+b+2,2)01} we have, 
\begin{itemize}
	\item either $-1-b+2a \pm 1= -b+2$, which imply either $2a-2b=0$ which is impossible by Lemma \ref{relations} or $4a-2b=0$ which is possible only for $n=7$ or $14$. However, direct SageMath computation for $n=7$ and $14$ show that such $\varphi$ does not exist.
	\item or $a-1+3b=-1-a \pm b$, i.e., $2a+4b=0$ which is possible only for $n=9$ or $2a+2b=0$, which is impossible by Lemma \ref{relations}. And finally direct SageMath computation for $n=9$ show that such $\varphi$ does not exist.
	\item or $-1-b+2a \pm 1=-3b$, i.e., $2a+2b=0$ or $2a+2b-2=0$, both of which are impossible by Lemma \ref{relations}. 
\end{itemize}
Hence $\varphi((a+b+2,2)) \neq (-b-a+2,2)$.

\textbf{Case B(3)(b):} If $\varphi((a+b+2,2))=(a-1-2b,1)$, then $\varphi((a+b+1,0)=(b-a-1,0)$ and $\varphi((2a+b,1)) \sim \varphi((a+b+2,2))$ imply 
\begin{equation}\label{(a+b+2,2)02}
\varphi((2a+b,1)) \in \lbrace (b-a-3,0), (1-b-2a \pm b,2) \rbrace.
\end{equation}
From the Equations \ref{(a+b+2,2)00} and \ref{(a+b+2,2)02} we have, 
\begin{itemize}
	\item either $b-a-3=-1+a \pm b$, i.e., $2a+2b=0$ or $2a-2b+2=0$, both of which are impossible by Lemma \ref{relations}.
	\item or $1-b-2a\pm b=-2a \pm 1$. These gives rise two four equations, out of which three are impossible by Lemma \ref{relations}, namely $2=0,~2b=0$ and $2a-2b=0$. The only possibility which remains is $1-b-2a+b=-2a+1$ and it is an identity.  
\end{itemize}
So assuming this to be the case, we have $\varphi((2b,0))=(-2,0)$, $\varphi((a+b+2,2))=(a-1-2b,1)$ and $\varphi((2a+b,1))=(-2a+1,2)$. 

Now consider the cycle $C_3: (2b,0) \sim (2a+b,1) \sim (a+b+2,2) \sim (1+a+3b,0) \sim (3a+1,1) \sim (3,2) \sim (2b,0)$. So $\varphi(C_3): (-2,0) \sim (-2a+1,2) \sim (a-1-2b,1) \sim \varphi((1+a+3b,0)) \sim \varphi((3a+1,1)) \sim \varphi((3,2)) \sim (-2,0)$. Now, $\varphi((2b,0))=(-2,0)$, $\varphi((2a+b,1))=(-2a+1,2)$ and $\varphi((3,2)) \sim (-2,0)$ imply $\varphi((3,2)) \in \lbrace (-3b,1), (-2a-1,2) \rbrace$. Again  $\varphi((2a+b,1))=(-2a+1,2)$, $\varphi((a+b+1,0))=(b-a-1,0)$ and $\varphi((1+a+3b,0)) \sim (a-1-2b,1)$ imply $\varphi((1+a+3b,0)) \in \lbrace (1-2a-2b,2), (b-a-3,0) \rbrace$.  Finally $\varphi((2b,0))=(-2,0)$ and $\varphi((3a+1,1)) \sim \varphi((3,2))$ imply 
	\begin{equation}\label{(3a+1,1)00}
	\varphi((3a+1,1)) \in \lbrace (-4,0), (-3a \pm b,2), (-2-2b,0), (-2b-a \pm 1,1) \rbrace.
	\end{equation}
\textbf{Case B(3)(b)(1):}	If $\varphi((1+a+3b,0))=(1-2a-2b,2)$ then $\varphi((3a+1,1)) \sim \varphi((1+a+3b,0))$ implies 
	\begin{equation}\label{(3a+1,1)01}
	\varphi((3a+1,1)) \in \lbrace (b-2a-2 \pm b,0), (a-2-2b \pm 1,1) \rbrace.
	\end{equation}
	From the Equations \ref{(3a+1,1)00} and \ref{(3a+1,1)01} we have, 
	\begin{itemize}
		\item either $b-2a-2 \pm b=-4$, i.e., $2a-2b=0$ or $2a-2b-2=0$, both of which are impossible by Lemma \ref{relations}.
		\item or $b-2a-2 \pm b= -2-2b$, i.e., $2a-2b=0$ or $2a-4b=0$,  both of which are impossible by Lemma \ref{relations}.
		\item or $a-2-2b \pm 1=-2b-a \pm 1$, i.e. $2a-2b=0$ or $2a=0$ or $4a-2b=0$. By Lemma \ref{relations}, the first two are impossible and the third one may hold only for $n=7$ or $14$. However, direct SageMath computation for $n=7$ and $14$ show that such $\varphi$ does not exist.
	\end{itemize}
	So we have $\varphi((1+a+3b,0))\neq (1-2b-2a,2)$.
	
\textbf{Case B(3)(b)(2):} If $\varphi((1+a+3b,0))=(b-a-3,0)$ then $\varphi((3a+1,1)) \sim \varphi((1+a+3b,0))$ implies 
	\begin{equation}\label{(3a+1,1)02}
	\varphi((3a+1,1)) \in \lbrace (a-1-3b \pm b,1), (b-1-3a \pm 1,2) \rbrace.
	\end{equation}
	From the Equations \ref{(3a+1,1)00} and \ref{(3a+1,1)02} we have, 
	\begin{itemize}
		\item either $a-1-3b \pm b= -2b -a \pm 1$, i.e., $2a=0$ or $2a-2b=0$ or $ 2a-2b-2=0$, all of which are impossible by Lemma \ref{relations}.
		\item or $1-b-3a \pm 1=-3a \pm b$. Out of the four relations that we get, three of them (namely, $2=0$, $2a-2b=0$ and $2b=0$) are invalid by Lemma \ref{relations} and the fourth is an identity, i.e., $1-b-3a - 1=-3a - b$.
	\end{itemize}
		
	So we have $\varphi((3a+1,1))=(-3a-b,2)$, $\varphi((3,2))=(-3b,1)$, $\varphi((a+1+3b,0))=(b-a-3,0)$.  Similarly we can show that $\varphi((3a-1,1))=(-3a+b,2)$ and  $\varphi((-a+1+3b,0))=(a+b-3,0)$.
	
Now $\varphi((2b,0))=(-2,0)$, $\varphi((3a+1,1))=(-3a-b,2)$, $\varphi((3a-1,1))=(-3a+b,2)$ and $\varphi((4b,0)) \sim \varphi((3,2))=(-3b,1) $ imply $\varphi((4b,0))=(-4,0)$.
	
Proceeding in this way, we can show that $\varphi((2kb,0))=(-2k,0) $, for all $k \in \Bbb Z$. So we have $\varphi((2,0))=(-2a,0)$, where $k=a$, which is a contradiction as we have shown earlier that $\varphi((2,0))=(2b,0)$ and $2b \neq -2a$. Therefore $\varphi((a+b+1,0)) \neq (b-a-1,0)$.

\textbf{Case B(1):}  When $\varphi((a+1,1))=(-a-b,2)$, $\varphi((1+a+b,0))=(2a+b,1)$, $\varphi((a+b,2))=(2b,0)$ then we have the Equation \ref{2a+2b+2=0,01}. As $2a+2b+2=0$, i.e., $a+b+1=-a-b-1$, then $\varphi((-1-a-b,0)) = (2a+b,1)$.  But from the Equations \ref{image-of-(-1-a-b,0)0} and \ref{image-of-(-1-a-b,0)1} we have $\varphi((-1-a-b,0)) \neq (2a+b,1)$, which is a contradiction. Hence $\varphi((1+a+b,0)) \neq (2a+b,1)$.

\textbf{Case B(2):}  Now when $\varphi((a+1,1))=(-a-b,2)$, $\varphi((a+b+1,0))=(a+b+1,0)$, $\varphi((a+b,2))=(a+1,1)$ then we have the Equation \ref{2a+2b+2=0,02}. Consider the cycle $C_1: (a+b+1,0) \sim (a+1,1) \sim (1+2b,2) \sim (2a,0) \sim (b+2,1) \sim (a+b,2) \sim (a+b+1,0)$. Then $\varphi(C_1): (a+b+1,0) \sim (-a-b,2) \sim \varphi((1+2b,2)) \sim \varphi((2a,0)) \sim \varphi((b+2,1)) \sim (a+1,1) \sim (a+b+1,0)$. Now $\varphi((a+b+1,0))=(a+b+1,0)=(-a-b-1,0)$, $\varphi((1,2))=(-b,1)$ and $\varphi((1+2b,2)) \sim (-a-b,2)$  imply $\varphi((1+2b,2)) \in \lbrace (-1-a+b,0), (-b-2,1) \rbrace$. Again, $\varphi((b,1))=(1,2)$, $\varphi((a+b+1,0))=(a+b+1,0)$ and $\varphi((b+2,1)) \sim \varphi((a+1,1))$ imply $\varphi((b+2,1)) \in \lbrace (b+a-1,0), (1+2b,2) \rbrace$. Also $\varphi((a+1,1))=(-a-b,2)$ and $\varphi((2a,0)) \sim \varphi((1+2b,2))$ imply 
\begin{equation}\label{(2a,0)00}
\varphi((2a,0)) \in \lbrace (-b-1+a \pm b,1), (-a-b+2,2), (-1-2a \pm 1,0), (-a-3b,2) \rbrace.
\end{equation}
\textbf{Case B(2)(a):}  Now if $\varphi((b+2,1))=(b+a-1,0)$ then $\varphi((a+b,2))=(a+1,1)$ and $\varphi((2a,0)) \sim \varphi((b+2,1))$ imply 
\begin{equation}\label{(2a,0)01}
\varphi((2a,0)) \in \lbrace (a+1-2b,1), (1+b-a \pm 1,2) \rbrace.
\end{equation}
From the Equations \ref{(2a,0)00} and \ref{(2a,0)01}, we have 
\begin{itemize}
\item either $1+b-a \pm 1=-a-b+2$, i.e., $2b=0$ or $2a-2b=0$, both of which are impossible by Lemma \ref{relations}.
\item or $1+b-a \pm 1= -a-3b$, i.e., $4b=0$ which is impossible or $4a+2b=0$, which, by Lemma \ref{relations}, holds only if $n=18$. However, direct SageMath computation for $n=18$ shows that such $\varphi$ does not exist.
\item or $a+1-2b= -b-1+a \pm b$, i.e., $2=0$ or $2a-2b=0$, both of which are impossible by Lemma \ref{relations}.  
\end{itemize}
Hence $\varphi((b+2,1)) \neq (b+a-1,0)$.

\textbf{Case B(2)(b):}  Now if $\varphi((b+2,1))=(1+2b,2)$ then $\varphi((a+b,2))=(a+1,1)$ and $\varphi((2a,0)) \sim \varphi((b+2,1))$ imply 
\begin{equation}\label{(2a,0)02}
\varphi((2a,0)) \in \lbrace (a+3,1), (b+2a \pm b,0) \rbrace.
\end{equation}
From the Equations \ref{(2a,0)00} and \ref{(2a,0)02}, we have 
\begin{itemize}
\item either $a+3=-b-1+a \pm b$, i.e., $4=0$ or $2a+4b=0$,  which, by Lemma \ref{relations}, can hold only if $n=9$. However, direct SageMath computation for $n=9$ shows that such $\varphi$ does not exist.
\item or $b+2a \pm b=-1-2a \pm 1$, i.e., $2(a+b+2)=0$ or $4a=0$ or $2a+4b=0$ or $4a+2b=0$. By Lemma \ref{relations}, the first two are impossible and the next two can hold only if $n=9$ or $18$. But those are also ruled out by SageMath computation. 
\end{itemize}
Hence $\varphi((b+2,1)) \neq (1+2b,2)$ and hence $\varphi((a+1,1))\neq (-a-b,2)$, i.e., {\bf Case B} can not hold.

{\bf Case C:} If $\varphi((a+1,1))=(-a+b,2)$ then $\varphi((1,2))=(-b,1)$ and $\varphi((a+b+1,0)) \sim \varphi((a+1,1))$ imply 
\begin{equation}\label{image-of-(1+a+b,0)13}
\varphi((1+a+b,0)) \in \lbrace (-b+2,1), (-1+a\pm b,0) \rbrace 
\end{equation}
From the Equations \ref{image-of-(1+a+b,0)10} and \ref{image-of-(1+a+b,0)13}, we have 
\begin{itemize}
\item either $-b+2=2a \pm b$, i.e., $2a-2b=0$ or $2(a+b-1)=0$, both of which are impossible by Lemma \ref{relations}.
\item or $-1+a\pm b=b-a\pm 1$, i.e., $2a=0$ or $2a-2b=0$ or $2(a-b-1)=0$ and all of them are ruled out by Lemma \ref{relations}. 
\item or $-1+a \pm b = b+a \pm 1$. This gives rise two four conditions, out of which three (namely, $2=0$, $2a+2b=0$ and $2b=0$) are ruled out by Lemma \ref{relations} and the fourth one is the identity $-1+a+b=b+a-1$.
\end{itemize}
 So we have $\varphi((a+b+1,0))=(a-1+b,0)$,  $\varphi(((a+b,2))=(a+1,1)$ and $\varphi((a+1,1))=(-a+b,2)$.
 
 Similarly we can show that $\varphi((-a+b+1,0))=(-a-1+b,0)$,  $\varphi(((a-b,2))=(a-1,1)$ and $\varphi((a-1,1))=(-a-b,2)$. 
 
 Now $\varphi((a+b,2))=(a+1,1)$, $\varphi(((a-b,2))=(a-1,1)$ and $\varphi((2,0)) \sim \varphi((b,1))=(1,2)$ imply $\varphi((2,0))=(2b,0)$.  

Consider the cycle $C_2: (a+b+1,0) \sim (a+1,1) \sim (1,2) \sim (2b,0) \sim (2a+b,1) \sim (a+b+2,2) \sim (a+b+1,0)$. So $\varphi(C_2): (a-1+b,0) \sim (-a+b,2) \sim (-b,1) \sim  \varphi((2b,0)) \sim \varphi((2a+b,1)) \sim \varphi((a+b+2,2)) \sim (a-1+b,0)$. Now,  $\varphi((0,0))=(0,0)$, $\varphi((a+1,0))=(-a+b,2)$ and $\varphi((2b,0)) \sim (-b,1)$ imply $\varphi((2b,0)) \in \lbrace (-a-b,2), (-2,0) \rbrace$. Again,  $\varphi(((a+b,2))=(a+1,1)$, $\varphi((a+1,1))=(-a+b,2)$ and $\varphi((a+b+2,2)) \sim  (a-1+b,0)$ imply $\varphi((a+b+2,2)) \in \lbrace (b-a+2,2) , (1-2b+a,1) \rbrace$. Also, $\varphi((1,2))=(-b,1)$ and $\varphi((2a+b,1)) \sim \varphi((2b,0)) $ imply 
\begin{equation}\label{(a+b+2,2)0}
\varphi((2a+b,1)) \in \lbrace (-1-a \pm b,0), (-b-2,1), (-3b,1), (-2a \pm 1,2) \rbrace.
\end{equation}
\textbf{Case C(1):}  If $\varphi((a+b+2,2))=(b-a+2,2)$ then $\varphi((a+b+1,0))=(a+b-1,0)$ and $\varphi((2a+b,1)) \sim \varphi((a+b+2,2))$ imply 
\begin{equation}\label{(a+b+2,2)1}
\varphi((2a+b,1)) \in \lbrace (a-1+3b,0), (1-b+2a \pm 1,1) \rbrace.
\end{equation}
From the Equations \ref{(a+b+2,2)0} and \ref{(a+b+2,2)1} we have
\begin{itemize}
\item either $1-b+2a \pm 1= -b-2$, i.e., $2a+2b=0$ or $4a+2b=0$. By Lemma \ref{relations}, the first is an impossibility and the second one can hold only if $n=18$. However, that is also ruled out by SageMath computation for $n=18$.
\item or $a-1+3b=-1-a \pm b$, i.e., $2a+2b=0$ or $2a+4b=0$. By Lemma \ref{relations}, the first is an impossibility and the second one can hold only if $n=9$. However, that is also ruled out by SageMath computation for $n=9$.
\item or $1-b+2a \pm 1=-3b$, i.e., $2a+2b=0$, which is impossible by the Lemma \ref{relations} but, 
\begin{equation}\label{2a+2b+2=0,03}
2a+2b+2=0~ may ~hold.
\end{equation} 
\end{itemize}

When $\varphi((a+b+1,0))=(a-1+b,0)$,  $\varphi(((a+b,2))=(a+1,1)$ and $\varphi((a+1,1))=(-a+b,2)$ then we have the Equation \ref{2a+2b+2=0,03}. As $2a+2b+2=0$, i.e., $a+b+1=-a-b-1$, then $\varphi((-1-a-b,0)) = (a-1+b,0)$.  But from the Equations \ref{image-of-(-1-a-b,0)0} and \ref{image-of-(-1-a-b,0)1}, we have $\varphi((-1-a-b,0)) \neq (a-1+b,0)$, which is a contradiction. Hence $\varphi((a+b+2,2)) \neq (b-a+2,2)$.

\textbf{Case C(2):} If $\varphi((a+b+2,2))=(a+1-2b,1)$, then $\varphi((a+b+1,0)=(a-1+b,0)$ and $\varphi((2a+b,1)) \sim \varphi((a+b+2,2))$ imply 
\begin{equation}\label{(a+b+2,2)2}
\varphi((2a+b,1)) \in \lbrace (b+a-3,0), (1+b-2a \pm b,2) \rbrace.
\end{equation}
From the Equations \ref{(a+b+2,2)0} and \ref{(a+b+2,2)2} we have, 
\begin{itemize}
\item either $b+a-3=-1-a \pm b$, i.e., $2a-2b=0$ or $2a+2b-2=0$, both of which are impossible by  Lemma \ref{relations}.
\item or $1+b-2a\pm b=-2a \pm 1$, This gives rise to four conditions, out of which three (namely $2b=0$, $2-0$ and $2a+2b=0$) are ruled out by Lemma \ref{relations} and the fourth one is the identity $1+b-2a-b=-2a+1$.
\end{itemize}
 So we have $\varphi((2a+b,1))=(-2a+1,2)$. Also previously, we had $\varphi((2b,0))=(-2,0)$ and $\varphi((a+b+2,2))=(1-2b+a,1)$.

Now consider the cycle $C_3: (2b,0) \sim (2a+b,1) \sim (a+b+2,2) \sim (1+a+3b,0) \sim (3a+1,1) \sim (3,2) \sim (2b,0)$. So $\varphi(C_3): (-2,0) \sim (-2a+1,2) \sim (1+a-2b,1) \sim \varphi((1+a+3b,0)) \sim \varphi((3a+1,1)) \sim \varphi((3,2)) \sim (-2,0)$. Now, $\varphi((2b,0))=(-2,0)$, $\varphi((2a+b,1))=(-2a+1,2)$ and $\varphi((3,2)) \sim (-2,0)$ imply $\varphi((3,2)) \in \lbrace (-3b,1), (-2a-1,2) \rbrace$. Again  $\varphi((2a+b,1))=(-2a+1,2)$, $\varphi((a+b+1,0))=(a+b-1,0)$ and $\varphi((1+a+3b,0)) \sim (1+a-2b,1)$ imply $\varphi((1+a+3b,0)) \in \lbrace (2b+1-2a,2), (a+b-3,0) \rbrace$.  Also, $\varphi((2b,0))=(-2,0)$ and $\varphi((3a+1,1)) \sim \varphi((3,2))$ imply 
\begin{equation}\label{(3a+1,1)0}
\varphi((3a+1,1)) \in \lbrace (-4,0), (-3a \pm b,2), (-2-2b,0), (-2b-a \pm 1,1) \rbrace.
\end{equation}

\textbf{Case C(2)(a):}  If $\varphi((1+a+3b,0))=(2b+1-2a,2)$, then $\varphi((3a+1,1)) \sim \varphi((1+a+3b,0))$ implies 
\begin{equation}\label{(3a+1,1)1}
\varphi((3a+1,1)) \in \lbrace (2a+b-2 \pm b,0), (2+a-2b \pm 1,1) \rbrace.
\end{equation}
From the Equations \ref{(3a+1,1)0} and \ref{(3a+1,1)1} we have
\begin{itemize}
\item either $2+a-2b \pm 1=-2b-a \pm 1$, i.e. $2a+2b=0$ or $2a=0$ (which are impossible by Lemma \ref{relations}) or $4a+2b=0$ which can hold only if $n=18$. But direct SageMath computation for $n=18$ ruled out this case.
\item or $2a+b-2 \pm b= -2-2b$, i.e., $2a+2b=0$ (impossible by Lemma \ref{relations}) or $2a+4b=0$, which  can hold only if $n=9$. But direct SageMath computation ruled out this possibility.
\item or $2a+b-2 \pm b=-4$, i.e., $2a+2b=0$, which is impossible by the Lemma \ref{relations} but 
\begin{equation}\label{2a+2b+2=0,04}
2a+2b+2=0 ~ may ~ hold.
\end{equation}	
\end{itemize}

Thus, if $\varphi((a+b+1,0))=(a-1+b,0)$,  $\varphi(((a+b,2))=(a+1,1)$ and $\varphi((a+1,1))=(-a+b,2)$ holds, then we have $2a+2b+2=0$. As $2a+2b+2=0$, i.e., $a+b+1=-a-b-1$, then $\varphi((-1-a-b,0)) = (a-1+b,0)$.  But from the Equations \ref{image-of-(-1-a-b,0)0} and \ref{image-of-(-1-a-b,0)1} we have $\varphi((-1-a-b,0)) \neq (a-1+b,0)$, which is a contradiction. Thus, Equation \ref{2a+2b+2=0,04} does not hold.  

So we have $\varphi((1+a+3b,0))\neq (2b+1-2a,2)$.

\textbf{Case C(2)(b):} If $\varphi((1+a+3b,0))=(a+b-3,0)$, then $\varphi((3a+1,1)) \sim \varphi((1+a+3b,0))$ implies 
\begin{equation}\label{(3a+1,1)2}
\varphi((3a+1,1)) \in \lbrace (1+a-3b \pm b,1), (b+1-3a \pm 1,2) \rbrace.
\end{equation}
From the Equations \ref{(3a+1,1)0} and \ref{(3a+1,1)2} we have, 
\begin{itemize}
\item either $1+a-3b \pm b= -2b -a \pm 1$, i.e., $2a=0$ or $2a+2b=0$ or $2a-2b=0$ or $ 2a-2b+2=0$, all of which are impossible by Lemma \ref{relations}.
\item or $b+1-3a \pm 1=-3a \pm b$. This gives rise to four conditions. Out of which three (namely, $2=0$,  $2a+2b=0$ and $2b=0$) are ruled out by Lemma \ref{relations} and fourth one is the identity $b+1-3a-1=-3a+b$. 
\end{itemize}
So we have $\varphi((3a+1,1))=(-3a+b,2)$, $\varphi((3,2))=(-3b,1)$, $\varphi((a+1+3b,0))=(a+b-3,0)$. 

Similarly we can show that $\varphi((3a-1,1))=(-3a-b,2)$ and $\varphi((-a+1+3b,0))=(-a+b-3,0)$.

Now $\varphi((2b,0))=(-2,0)$, $\varphi((3a+1,1))=(-3a+b,2)$, $\varphi((3a-1,1))=(-3a-b,2)$ and $\varphi((4b,0)) \sim \varphi((3,2))=(-3b,1) $ imply $\varphi((4b,0))=(-4,0)$.

Proceeding this way, we can show that $\varphi((2kb,0))=(-2k,0) $, for all $k \in \Bbb Z$. So we have $\varphi((2,0))=(-2a,0)$, where $k=a$, which is a contradiction as we have shown that $\varphi((2,0))=(2b,0)$ and $2b \neq -2a$. Therefore we have $\varphi((a+1,1))\neq (-a+b,2)$ and {\bf Case C} can not hold.

As none of the {\bf Cases A, B} and {\bf C} hold, the assumption that $\varphi((1,2)) = (-b,1)$ is wrong. Hence the lemma follows. \qed

{\theorem \label{phi(1,2)=(b,1)} If $\varphi \in G$ and $\varphi((0,0))=(0,0)$, $\varphi((b,1))=(1,2)$, $\varphi((1,2))=(b,1)$ then $n=7$ or $14$. }\\
\pf  Consider the cycle $C_0: (0,0)\sim (b,1)\sim (a+b,2)\sim (a+b+1,0) \sim (a+1,1)\sim (1,2)\sim (0,0)$. Then $\varphi(C_0): (0,0)\sim (1,2)\sim \varphi((a+b,2)) \sim \varphi((a+b+1,0)) \sim \varphi((a+1,1)) \sim (b,1) \sim(0,0) \rbrace$. $\varphi((0,0))=(0,0)$ and $\varphi((a+b,2)) \sim (1,2)$ imply $\varphi((a+b,2)) \in \lbrace (2b,0), (a \pm 1,1)$. $\varphi((0,0))=(0,0)$ and $\varphi((a+1,1)) \sim (b,1)$ imply $\varphi((a+1,1)) \in \lbrace (2,0), (a \pm b,2) \rbrace$. Now $\varphi((b,1))=(1,2)$ and $\varphi((a+b+1,0)) \sim \varphi((a+b,2))$ imply 
\begin{equation}\label{image-of-(1+a+b,0)20}
\varphi((1+a+b,0)) \in \lbrace (2a \pm b,1), (3,2), (1+ 2b,2), (b +a \pm 1, 0),(1-2b,2), (b-a \pm 1,0)  \rbrace.
\end{equation}
Depending upon the value of $\varphi((a+1,1))$, we split it into three cases {\bf a, b} and {\bf c}.

\textbf{Case a:}  If $\varphi((a+1,1))=(2,0)$ then $\varphi((1+a+b,0)) \sim \varphi((a+1,1))$ and $\varphi((1,2))=(-b,1))$ imply
\begin{equation}\label{image-of-(1+a+b,0)21}
\varphi((1+a+b,0)) \in \lbrace (3b,1), (2a \pm 1,2) \rbrace 
\end{equation}
From the Equation \ref{image-of-(1+a+b,0)20} and \ref{image-of-(1+a+b,0)21} we have,
\begin{itemize}
	\item either $3b=2a \pm b$, i.e., $2a-2b=0$ or $4b-2a=0$, which are impossible by Lemma \ref{relations}.
	\item or $2a \pm 1=1+2b$, i.e, $2a-2b=0$ or $2a-2b-2=0$, which are impossible by Lemma \ref{relations}.
	\item or $2a \pm 1=1-2b$, i.e, $2a+2b=0$ or $2a+2b-2=0$, which are impossible by Lemma \ref{relations}.
	\item or $2a \pm 1=3$, i.e., $2a=2$, i.e., $2a-2b=0$ or  $2a=4$, i.e., $4a-2b=0$. By Lemma \ref{relations}, $2a-2b=0$ can not hold and $4a-2b=0$ can hold only if $n=7, 14$. However direct SageMath computation for $n=7,14$ shows that such $\varphi$ does not exist.
\end{itemize}  Therefore we have $\varphi((a+1,1)) \neq (2,0)$ and hence {\bf Case a} can not hold. 

\textbf{Case b:} If $\varphi((a+1,1))=(a+b,2)$ then  $\varphi((1+a+b,0)) \sim \varphi((a+1,1))$ and $\varphi((1,2))=(b,1)$  imply
\begin{equation}\label{image-of-(1+a+b,0)22}
\varphi((1+a+b,0)) \in \lbrace (b+2,1), (1+a\pm b,0) \rbrace 
\end{equation}
From the Equation \ref{image-of-(1+a+b,0)20} and \ref{image-of-(1+a+b,0)22} we have,
\begin{itemize}
	\item either $b+2=2a \pm b$, i.e., $2a=2$, i.e., $2a-2b=0$ or $2a-2b-2=0$, which are impossible.
	\item or $1+a \pm b=b-a \pm 1$, i.e., $2a=0$ or $2a+2=0$, i.e., $2a+2b=0$ or $2a-2b=0$ or $2a-2b+2=0$, which are impossible.
	\item or $1+a \pm b=b+a \pm 1$. This gives rise four equations out of which three are impossible by Lemma \ref{relations}, namely $2=0$, $2b=0$, $2b-2=0$, i.e., $2a-2b=0$. The only possiblity which remains is $1+a+b=b+a+1$, which is an identity. So asssuming this to be the case we have $\varphi((a+1,1))=(a+b,2)$, $\varphi((a+b+1,0))=(a+b+1,0)$ and $\varphi((a+b,2))=(a+1,1)$.
\end{itemize} 

\textbf{Case c:}  If $\varphi((a+1,1))=(a-b,2)$ then  $\varphi((1+a+b,0)) \sim \varphi((a+1,1))$ and $\varphi((1,2))=(b,1)$  imply
\begin{equation}\label{image-of-(1+a+b,0)23}
\varphi((1+a+b,0)) \in \lbrace (b-2,1), (1-a\pm b,0) \rbrace 
\end{equation}
From the Equation \ref{image-of-(1+a+b,0)20} and \ref{image-of-(1+a+b,0)23} we have,
\begin{itemize}
	\item either $b-2=2a \pm b$, i.e., $2a+2=0$, i.e., $2a+2b=0$ or $2a-2b+2=0$, which are impossible.
	\item or $1-a \pm b=b+a \pm 1$, i.e., $2=0$ or $2a-2=0$, i.e., $2a-2b=0$ or $2a+2b=0$ or $2a+2b-2=0$, which are impossible.
	\item or $1-a \pm b=b-a \pm 1$,This gives rise four equations out of which three are impossible by Lemma \ref{relations}, namely $2=0$, $2b=0$, $2b-2=0$, i.e., $2a-2b=0$. The only possiblity which remains is $1-a+b=b-a+1$, which is an identity. So asssuming this to be the case we have $\varphi((a+1,1))=(a-b,2)$, $\varphi((a+b+1,0))=(b-a+1,0)$ and $\varphi((a+b,2))=(a-1,1)$.
\end{itemize} 
Combining the feasible cases in Case {\bf b} and {\bf Case b}, we have  $\varphi((a+1,1)) \in \lbrace (a \pm b,2) \rbrace$, $\varphi((a+b+1,0))\in \lbrace (b\pm a +1,0 ) \rbrace$ and $\varphi((a+b,2)) \in \lbrace (a \pm 1,1) \rbrace$.

Now consider the cycle $C'_0: (0,0)\sim (1,2)\sim (a-b,2)\sim (b-a+1,0) \sim (a-1,1)\sim (1,2)\sim (0,0)$. Then $\varphi(C'_0): (0,0)\sim (1,2)\sim \varphi((a-b,2)) \sim \varphi((b-a+1,0)) \sim \varphi((a-1,1)) \sim (b,1) \sim(0,0) \rbrace$. $\varphi((0,0))=(0,0)$ and $\varphi((a-b,2)) \sim (1,2)$ imply $\varphi((a-b,2)) \in \lbrace (2b,0), (a \pm 1,1)$. $\varphi((0,0))=(0,0)$ and $\varphi((a-1,1)) \sim (b,1)$ imply $\varphi((a-1,1)) \in \lbrace (2,0), (a \pm b,2) \rbrace$. Now $\varphi((b,1))=(1,2)$ and $\varphi((b-a+1,0)) \sim \varphi((a-b,2))$ imply 
\begin{equation}\label{image-of-(b-a+1,0)20}
\varphi((b-a+1,0)) \in \lbrace (2a \pm b,1), (3,2), (1+ 2b,2), (b +a \pm 1, 0),(1-2b,2), (b-a \pm 1,0)  \rbrace.
\end{equation}
\textbf{Case a$^\prime$:} If $\varphi((a-1,1))=(2,0)$ then $\varphi((b-a+1,0)) \sim \varphi((a-1,1))$ and $\varphi((1,2))=(-b,1))$ imply
\begin{equation}\label{image-of-(b-a+1,0)21}
\varphi((b-a+1,0)) \in \lbrace (3b,1), (2a \pm 1,2) \rbrace 
\end{equation}
From the Equation \ref{image-of-(b-a+1,0)20} and \ref{image-of-(b-a+1,0)21} we have,
\begin{itemize}
	\item either $3b=2a \pm b$, i.e., $2a-2b=0$ or $4b-2a=0$, which are impossible by the Lemma \ref{relations}.
	\item or $2a \pm 1=1+2b$, i.e, $2a-2b=0$ or $2a-2b-2=0$, which are impossible by the Lemma \ref{relations}.
	\item or $2a \pm 1=1-2b$, i.e, $2a+2b=0$ or $2a+2b-2=0$, which are impossible by the Lemma \ref{relations}.
	\item or $2a \pm 1=3$, i.e., $2a=2$, i.e., $2a-2b=0$ or  $2a=4$, i.e., $4a-2b=0$. By Lemma \ref{relations} $2a-2b=0$ can not hold and $4a-2b=0$ can hold only if $n=7, 14$. Therefore we have $\varphi((a-1,1))=(2,0)$ only if $n=7,14$. 
\end{itemize} 

\textbf{Case b$^\prime$:} If $\varphi((a-1,1))=(a+b,2)$ then  $\varphi((b-a+1,0)) \sim \varphi((a-1,1))$ and $\varphi((1,2))=(b,1)$  imply
\begin{equation}\label{image-of-(b-a+1,0)22}
\varphi((b-a+1,0)) \in \lbrace (b+2,1), (1+a\pm b,0) \rbrace 
\end{equation}
From the Equation \ref{image-of-(b-a+1,0)20} and \ref{image-of-(b-a+1,0)22} we have,
\begin{itemize}
	\item either $b+2=2a \pm b$, i.e., $2a=2$, i.e., $2a-2b=0$ or $2a-2b-2=0$, which are impossible.
	\item or $1+a \pm b=b-a \pm 1$, i.e., $2a=0$ or $2a+2=0$, i.e., $2a+2b=0$ or $2a-2b=0$ or $2a-2b+2=0$, which are impossible.
	\item or $1+a \pm b=b+a \pm 1$. This gives rise four equations out of which three are impossible by Lemma \ref{relations}, namely $2=0$, $2b=0$, $2b-2=0$, i.e., $2a-2b=0$. The only possiblity which remains is $1+a+b=b+a+1$, which is an identity. So asssuming this to be the case we have $\varphi((a-1,1))=(a+b,2)$, $\varphi((b-a+1,0))=(a+b+1,0)$ and $\varphi((a-b,2))=(a+1,1)$.
\end{itemize} 

\textbf{Case c$^\prime$:}  If $\varphi((a-1,1))=(a-b,2)$ then  $\varphi((b-a+1,0)) \sim \varphi((a-1,1))$ and $\varphi((1,2))=(b,1)$  imply
\begin{equation}\label{image-of-(b-a+1,0)23}
\varphi((b-a+1,0)) \in \lbrace (b-2,1), (1-a\pm b,0) \rbrace 
\end{equation}
From the Equation \ref{image-of-(b-a+1,0)20} and \ref{image-of-(b-a+1,0)23} we have,
\begin{itemize}
	\item either $b-2=2a \pm b$, i.e., $2a+2=0$, i.e., $2a+2b=0$ or $2a-2b+2=0$, which are impossible.
	\item or $1-a \pm b=b+a \pm 1$, i.e., $2=0$ or $2a-2=0$, i.e., $2a-2b=0$ or $2a+2b=0$ or $2a+2b-2=0$, which are impossible.
	\item or $1-a \pm b=b-a \pm 1$, this gives rise four equations out of which three are impossible by Lemma \ref{relations}, namely $2=0$, $2b=0$, $2b-2=0$, i.e., $2a-2b=0$. The only possiblity which remains is $1-a+b=b-a+1$, which is an identity. So asssuming this to be the case we have $\varphi((a-1,1))=(a-b,2)$, $\varphi((b-a+1,0))=(b-a+1,0)$ and $\varphi((a-b,2))=(a-1,1)$.
\end{itemize} 
Combining these three cases we have  $\varphi((a-1,1)) \in \lbrace (a \pm b,2) \rbrace$, $\varphi((b-a+1,0)) \in \lbrace (b\pm a +1,0 ) \rbrace$,  $\varphi((a-b,2)) \in \lbrace (a \pm 1,1) \rbrace$.    $\varphi((a-1,1))=(2,0) $, $\varphi((b-a+1,0))=(3,2)$, $\varphi((a-b,2))=(2b,0)$ only if $n=7,14$. 

Depending upon the values of $\varphi((a+1,1))$ and $\varphi((a-1,1))$ (that we get from Cases {\bf a, b, c} and {\bf a$^\prime$, b$^\prime$, c$^\prime$}), we have four different cases namely \textbf{Case A:} $\varphi((a+1,1))=(a+b,2)$ and $\varphi((a-1,1))=(2,0)$ (only if $n=7,14$), \textbf{Case B:} $\varphi((a+1,1))=(a-b,2)$ and $\varphi((a-1,1))=(2,0)$ (only if $n=7,14$), \textbf{Case C:} $\varphi((a+1,1))=(a-b,2)$ and $\varphi((a-1,1))=(a+b,2)$ and \textbf{Case D:} $\varphi((a+1,1))=(a+b,2)$ and $\varphi((a-1,1))=(a-b,2)$. However, before resolving these four cases, we prove a claim which will be crucial in the following proof.

\textbf{Claim:}  $\varphi((-1-a-b,0)) \in \lbrace (-1+a \pm b,0), (-b \pm 2,1), (-1-a \pm b,0)$,

\hspace{5.5cm} $(-2a \pm 1,2), (-3b,1), (-b+a \pm 1,0,0),(-1 \pm 2b,2) $,

\hspace{5.5cm} $ (-b-a \pm 1,0), (-2a \pm b,1),  (-3,2) \rbrace.$

\textbf{Proof of Claim:} $\varphi((0,0))=(0,0)$, $\varphi((b,1))=(1,2)$, $\varphi((1,2))=(b,1)$, $\varphi((-b,1)) \sim \varphi((0,0))=(0,0)$ and  imply $\varphi((-1,2)) \sim \varphi((0,0))=(0,0)$ imply $\varphi((-b,1)),~\varphi((-1,2))\in \lbrace (-b,1), (-1,2) \rbrace $.

\textbf{Case 1:} Let $\varphi((-b,1))=(-b,1)$ and $\varphi((-1,2))=(-1,2)$. Consider the cycle $C''_0: (0,0)\sim (-b,1)\sim (-a+-,2)\sim (-a-b-1,0) \sim (-a-1,1)\sim (-1,2)\sim (0,0)$. Then $\varphi(C''_0): (0,0)\sim (-b,1)\sim \varphi((-a-b,2)) \sim \varphi((-a-b-1,0)) \sim \varphi((-a-1,1)) \sim (-1,2) \sim(0,0) \rbrace$.  $\varphi((-a-b,2)) \sim (-b,1)$ and $\varphi((0,0))=(0,0)$ imply $\varphi((-a-b,2)) \in \lbrace (-a \pm b,2), (-2,0) \rbrace$.  $\varphi((-a-1,1)) \sim (-1,2)$ and $\varphi((0,0))=(0,0)$ imply $\varphi((-a-1,1)) \in \lbrace (-a \pm 1,1), (-2b,0) \rbrace$. $\varphi((-a-b-1,0)) \sim \varphi ((-a-b,2))$ and $\varphi((-b,1))=(-b,1)$ imply 
\begin{equation}\label{(-a-b-1,0)20}
\varphi((-a-b-1,0)) \in \lbrace (-1+a \pm b,0), (-b+2,1), (-1-a \pm b,0), (-b-2,1), (-2a \pm 1,2), (-3b,1) \rbrace.
\end{equation}
\textbf{Case 2:} Let $\varphi((-b,1))=(-1,2)$ and $\varphi((-1,2))=(-b,1)$. Again consider the cycle $C''_0: (0,0)\sim (-b,1)\sim (-a+-,2)\sim (-a-b-1,0) \sim (-a-1,1)\sim (-1,2)\sim (0,0)$. Then $\varphi(C''_0): (0,0)\sim (-1,2)\sim \varphi((-a-b,2)) \sim \varphi((-a-b-1,0)) \sim \varphi((-a-1,1)) \sim (-1,2) \sim(0,0) \rbrace$.  $\varphi((-a-b,2)) \sim (-1,2)$ and $\varphi((0,0))=(0,0)$ imply $\varphi((-a-b,2)) \in \lbrace (-a \pm 1,1), (-2b,0) \rbrace$. $\varphi((-a-1,1)) \sim (-b,1)$ and $\varphi((0,0))=(0,0)$ imply $\varphi((-a-1,1)) \in \lbrace (-a \pm b,2), (-2,0) \rbrace$.  $\varphi((-a-b-1,0)) \sim \varphi ((-a-b,2))$ and $\varphi((-b,1))=(-1,2)$ imply 
\begin{equation}\label{(-a-b-1,0)21}
\varphi((-a-b-1,0)) \in \lbrace (-b+a\pm 1,0), (-1+2b,2), (-b-a\pm 1,0), (-1-2b,2), (-2a \pm b,1), (-3,2) \rbrace.
\end{equation}

\textbf{Case A:} Let $\varphi((a+1,1))=(a+b,2)$, $\varphi((a+b+1,0))= (b+a +1,0 ) $, $\varphi((a+b,2))= (a+1,1)$, $\varphi((a-1,1))=(2,0)$, $\varphi((b-a+1,0))=(3,2)$ and $\varphi((a-b,2))=(2b,0)$. This map exists only if $n=7,14$. This can also be checked with direct Sagemath computation.

\textbf{Case B:} Let $\varphi((a+1,1))=(a-b,2)$, $\varphi((a+b+1,0))= (b-a +1,0 ) $, $\varphi((a+b,2))= (a-1,1)$, $\varphi((a-1,1))=(2,0)$, $\varphi((b-a+1,0))=(3,2)$ and $\varphi((a-b,2))=(2b,0)$. Such $\varphi$ can hold only if $n=7,14$. However direct SageMath computation for $n=7,14$ shows that such $\varphi$ does not exist, hence this map does not exist for all $n$.

\textbf{Case C:} Let $\varphi((a+1,1))=(a-b,2)$, $\varphi((a+b+1,0))= (b-a +1,0 ) $, $\varphi((a+b,2))= (a-1,1)$, $\varphi((a-1,1))=(a+b,2)$, $\varphi((b-a+1,0))= (b+a +1,0 ) $, $\varphi((a-b,2))= (a+1,1)$. Now  $\varphi((a+b,2))= (a-1,1)$, $\varphi((a-b,2))= (a+1,1)$ and $\varphi((2,0)) \sim \varphi((b,1))=(1,2)$ imply $\varphi((2,0))=(2b,0)$. 

Consider the cycle $C_1: (2,0) \sim (b,1) \sim (a+b,2) \sim (a+b+1,0) \sim (2b+a+1,1) \sim (2a+1,2) \sim (2,0)$. then $\varphi(C_1): (2b,0) \sim (1,2) \sim (a-1,1) \sim (b-a+1,0) \sim \varphi((2b+a+1,1)) \sim \varphi((2a+1,2)) \sim (2b,0)$. $\varphi((2b+a+1,1)) \sim (b-a+1,0)$, $\varphi((a+1,1))=(a-b,2)$ and $\varphi((a+b,2))=(a-1,1)$ imply $\varphi((2b+a+1,1)) \in \lbrace (a-1+2b,1), (2-b+a,2)  \rbrace$. Now $\varphi((2a+1,2)) \sim (2b,0)$ and $\varphi((b,1))=(1,2)$ imply 
\begin{equation}\label{(2a+1),20}
\varphi((2a+1,2)) \in \lbrace (2a \pm b,1), (3,2) \rbrace.
\end{equation}
\textbf{Case C(a):}  If $\varphi((2b+a+1,1))=(a-1+2b,1)$ then $\varphi((2b+a+1,1)) \sim \varphi((2a+1,2))$ and $\varphi((a+b+1,0))=(b-a+1,0)$ imply 
\begin{equation}\label{(2a+1),21}
\varphi((2a+1,2)) \in \lbrace (1-b+2a \pm b,2), (b-a+3,0) \rbrace.
\end{equation}
From the Equation \ref{(2a+1),20} and \ref{(2a+1),21} we have, $1-b+2a \pm b=3$, i.e., $2a-2b-2=0$ or $2a-2=0$, i.e., $2a-2b=0$, which are impossible by the Lemma \ref{relations}, hence $\varphi((2b+a+1,1))\neq (a-1+2b,1)$.\\
\textbf{Case C(b):}  If $\varphi((2b+a+1,1))= (2-b+a,2)$ then $\varphi((2b+a+1,1)) \sim \varphi((2a+1,2))$ and $\varphi((a+b+1,0))=(b-a+1,0)$ imply 
\begin{equation}\label{(2a+1),22}
\varphi((2a+1,2)) \in \lbrace (2a-1+b \pm 1,1), (3b-a+1,0) \rbrace.
\end{equation}
From the Equation \ref{(2a+1),20} and \ref{(2a+1),22} we have, $2a-1+b \pm 1=2a \pm b$, This gives rise four equations out of which three are impossible by Lemma \ref{relations}, namely $2=0$, $2b=0$, $2b-2=0$, i.e., $2a-2b=0$. The only possiblity which remains is $2a+b=2a+b$, which is an identity. So asssuming this to be the case we have $\varphi((2a+1,2))=(2a+b,1)$, $\varphi((2b+a+1,1))= (2-b+a,2)$.

Now consider the cycle $C'_1: (2,0) \sim (b,1) \sim (a+b,2) \sim (1+a-b,0) \sim (2b+1-a,1) \sim (2a-1,2) \sim (2,0)$. Then $\varphi(C'_1): (2b,0) \sim (1,2) \sim (a-1,1) \sim \varphi((1+a-b,0)) \sim \varphi((2b+1-a,1)) \sim \varphi((2a-1,2)) \sim (2b,0)$. $\varphi((2a-1,2)) \sim (2b,0)$, $\varphi((b,1))=(1,2)$ and $\varphi((2a+1,2))=(2a+b,1)$ imply $\varphi((2a-1,2)) \in \lbrace (2a-b,1), (3,2) \rbrace.$ Now  $ \varphi((1+a-b,0)) \sim(a-1,1)$, $\varphi((a+b+1,0))=(b-a+1,0)$ and $\varphi((b,1))=(1,2)$ imply $\varphi((1+a-b,0)) \in \lbrace (b-a-1,0), (1-2b,2) \rbrace$. $\varphi((2b+1-a,1)) \sim \varphi((1+a-b,0))$ and $\varphi((a+b,2))=(a-1,1)$ imply
\begin{equation}\label{(2b+1-a,1)20}
\varphi((2b+1-a,1)) \in \lbrace  (a-1-2b,1), (1-b-a \pm 1,2), (a-3,1), (b-2a \pm b,0) \rbrace.
\end{equation}
\textbf{Case C(b)(1):}  If $\varphi((2a-1,2))=(3,2)$ then $\varphi((2b+1-a,1)) \sim \varphi((2a-1,2))$ and $\varphi((2,0))=(2b,0)$ imply
\begin{equation}\label{(2b+1-a,1)21}
\varphi((2b+1-a,1)) \in \lbrace  (4b,0), (3a \pm 1,1)\rbrace.
\end{equation}
From the Equations \ref{(2b+1-a,1)20} and \ref{(2b+1-a,1)21} we have,
\begin{itemize}
	\item either $4b=b-2a \pm b$, i.e., $2a+2b=0$ or $2a+4b=0$, which are impossible by the Lemma \ref{relations}.
	\item $3a \pm 1=a-3$, i.e., $2a+4=0$, i.e., $4a+2b=0$ or $2a+2=0$, i.e., $2a+2b=0$. By Lemma \ref{relations} $2a+2b=0$ is impossible and $4a+2b=0$ holds only if $n=18$. But this is also ruled out by SageMath computation for $n=18$.
	\item or $3a \pm 1=a-1-2b$, i.e., $2a+2b=0$, which is impossible by the Lemma \ref{relations}, but $2a+2b+2=0$ may hold. We have this Equation when $\varphi((a+b+1,0))= (b-a +1,0 ) $. So $2a+2b+2=0$ imply $a+b+1=-a-b-1$, hence $\varphi((-a-b-1,0))= (b-a +1,0 ) $. But from the Equations \ref{(-a-b-1,0)20} and \ref{(-a-b-1,0)21} we have $\varphi((-a-b-1,0)) \neq (b-a +1,0 ) $, which is a contradiction. Hence $\varphi((2a-1,2))\neq (3,2)$.
\end{itemize}  

\textbf{Case C(b)(2):}  If $\varphi((2a-1,2))=(2a-b,1)$ then $\varphi((2b+1-a,1)) \sim \varphi((2a-1,2))$ and $\varphi((2,0))=(2b,0)$ imply
\begin{equation}\label{(2b+1-a,1)22}
\varphi((2b+1-a,1)) \in \lbrace  (2b-2,0), (2-a \pm b,2)\rbrace.
\end{equation}
From the Equations \ref{(2b+1-a,1)20} and \ref{(2b+1-a,1)22} we have, \begin{itemize}
	\item either $2b-2=b-2a \pm b$, i.e., $2a-2=0$, i.e., $2a-2b=0$ or $2b+2a-2=0$, which are impossible by the Lemma \ref{relations}.\\
	\item or $2-a \pm b= 1-b -a \pm 1$, This gives rise four equations out of which three are impossible by Lemma \ref{relations}, namely $2=0$, $2b=0$, $2b+2=0$, i.e., $2a+2b=0$. The only possiblity which remains is $2-a-b=2-b-a$, which is an identity. So asssuming this to be the case we have $\varphi((2b+1-a,0))=(2-b-a,0)$, $\varphi((2a-1,2))=(2a-b,1)$ and $\varphi((1+a-b,0))=(b-a-1,0)$.
\end{itemize}

Now $\varphi((3b,1)) \sim \varphi((2,0))=(2b,0) $, $\varphi((2a+1,2))=(2a+b,1)$ and $\varphi((2a-1,2))=(2a-b,1)$ imply  $\varphi((3b,1)) =(3,2)$.

Consider the cycle $C_2: (2,0) \sim ( 3b,1) \sim (3a+b,2) \sim (3+a+b,0) \sim (2b+a+1,1) \sim (2a+1,2) \sim (2,0)$. Then $\varphi(C_2): (2b,0) \sim (3,2) \sim \varphi((3a+b,2) \sim \varphi((3+a+b,0)) \sim (a-b+2,2) \sim (2a+b,1) \sim (2b,0)$. $\varphi((3a+b,2)) \sim (3,2)$  and $\varphi((2,0))=(2b,0)$ imply $\varphi((3a+b,2)) \in \lbrace (4b,0), (3a \pm 1,1) \rbrace$. $\varphi((3+a+b,0)) \sim (a-b+2,2) $, $\varphi((a+b+1,0))=(b-a+1,0)$ and $\varphi((2a+1,2))=(2a+b,1)$ imply
\begin{equation}\label{(3+a+b,0)20}
\varphi((3+a+b,0)) \in \lbrace (2a+b-2,1), (3b-a+1,0) \rbrace.
\end{equation}
\textbf{Case C(b)(2)(a):}  If $\varphi((3a+b,2))=(4b,0)$ then $\varphi((3+a+b,0)) \sim \varphi((3a+b,2))$ and $\varphi((3b,1)) =(3,2)$ imply 
\begin{equation}\label{(3+a+b,0)21}
\varphi((3+a+b,0)) \in \lbrace (4a \pm b,1), (5,2) \rbrace.
\end{equation}
From the Equations \ref{(3+a+b,0)20} and \ref{(3+a+b,0)21} we have, $4a \pm b= 2a +b-2$, i.e., $2a+2=0$, i.e., $2a+2b=0$, which is impossible by the Lemma \ref{relations} but $2a+2b+2=0$ may hold. We have this Equation when $\varphi((a+b+1,0))= (b-a +1,0 ) $. So $2a+2b+2=0$ imply $a+b+1=-a-b-1$, hence $\varphi((-a-b-1,0))= (b-a +1,0 ) $. But from the Equations \ref{(-a-b-1,0)20} and \ref{(-a-b-1,0)21} we have $\varphi((-a-b-1,0)) \neq (b-a +1,0 ) $, which is a contradiction. Hence $\varphi((3a+b,2))\neq (4b,0)$.

\textbf{Case C(b)(2)(b):}  If $\varphi((3a+b,2))=(3a+1,1)$  then $\varphi((3+a+b,0)) \sim \varphi((3a+b,2))$ and $\varphi((3b,1)) =(3,2)$ imply 
\begin{equation}\label{(3+a+b,0)22}
\varphi((3+a+b,0)) \in \lbrace (3b+a \pm 1,0), (3+2b,2) \rbrace.
\end{equation}
From the Equations \ref{(3+a+b,0)20} and \ref{(3+a+b,0)22} we have, $3b +a \pm 1=3b-a+1$, i.e., $2a=0$ or $2a-2=0$, i.e., $2a-2b=0$, which is impossible by the Lemma \ref{relations}. Hence $\varphi((3a+b,2))\neq (3a+1,1)$.

\textbf{Case C(b)(2)(c):}  If $\varphi((3a+b,2))=(3a-1,1)$  then $\varphi((3+a+b,0)) \sim \varphi((3a+b,2))$ and $\varphi((3b,1)) =(3,2)$ imply 
\begin{equation}\label{(3+a+b,0)23}
\varphi((3+a+b,0)) \in \lbrace (3b-a \pm 1,0), (3-2b,2) \rbrace.
\end{equation}
From the Equations \ref{(3+a+b,0)20} and \ref{(3+a+b,0)23} we have, $3b-a \pm 1=3b-a+1$, this gives rise to two equations, out of one is impossible, namely $2=0$. The only possibility which remains is $3b-a+1=3b-a+1$, which is an identity. So assuming this to be the case, we have $\varphi((3a+b,2))=(3a-1,1)$ and $\varphi((3+a+b,0))=(3b-a+1,0)$.

Now consider the cycle $C'_2: (2,0) \sim (3b,1) \sim (3a-b,2) \sim (3-a-b,0) \sim (2b-a-1,1) \sim (2a-1,2) \sim (2,0)$. $\varphi(C'_2): (2b,0) \sim (3,2) \sim \varphi((3a-b,2)) \sim \varphi((3-a-b,0)) \sim \varphi((2b-a-1,1)) \sim (2a-b,1) \sim (2b,0)$. $\varphi((3a-b,2)) \sim (3,2)$, $\varphi((a-b,2))=(a+1,1)$ and $\varphi((2,0))=(2b,0)$ imply $\varphi((3a-b,2)) \in \lbrace (3a+1,1), (4b,0) \rbrace$. $\varphi((2b-a-1,1)) \sim (2a-b,1)$ and $\varphi((2,0))=(2b,0)$ imply $\varphi((2b-a-1,1)) \in \lbrace (2-a \pm b,2), (2b-2,0) \rbrace$. Now $\varphi((3-a-b,0)) \sim \varphi((3a-b,2)) $ and $\varphi((3b,1))=(5,2)$ imply
\begin{equation}\label{(3-a-b,0),20}
\varphi((3-a-b,0)) \in \lbrace (4a \pm b,1), (5,2), (3+2b,2), (3b+a \pm 1,0) \rbrace.
\end{equation}
\textbf{Case C(b)(2)(d):}  If $\varphi((2b-a-1,1))=(2b-2,0)$ then $\varphi((3-a-b,0)) \sim \varphi((2b-a-1,1))$ and $\varphi((2a-1,2))=(2a-b,1)$ imply
\begin{equation}\label{(3-a-b,0),21}
\varphi((3-a-b,0)) \in \lbrace (2a-3b,1), (2-2a \pm 1,2) \rbrace.
\end{equation}
From the Equations \ref{(3-a-b,0),20} and \ref{(3-a-b,0),21} we have,
\begin{itemize}
	\item either $2a-3b=4a \pm b$, i.e., $2a+2b=0$ or $2a+4b=0$. By Lemma \ref{relations} $2a+2b=0$ can not hold and $2a+4b=0$ can hold only if $n=9$. However, direct SageMath computation for $n=9$ shows that such $\varphi$ does not exist.
	\item or $2-2a \pm 1=5$, i.e., $2a+2=0$, i.e., $2a+2b=0$ or $2a+4=0$, i.e., $4a+2b=0$. By Lemma \ref{relations} $2a+2b=0$ can not hold and $4a+2b=0$ can hold only if $n=18$. However, direct SageMath computation for $n=18$ shows that such $\varphi$ does not exist.
	\item or $2-2a \pm 1=3+2b$, i.e., $2a+2b=0$, which is impossible by the Lemma \ref{relations} but $2a+2b+2=0$ may hold. We have this Equation when $\varphi((a+b+1,0))= (b-a +1,0 ) $. So $2a+2b+2=0$ imply $a+b+1=-a-b-1$, hence $\varphi((-a-b-1,0))= (b-a +1,0 ) $. But from the Equations \ref{(-a-b-1,0)20} and \ref{(-a-b-1,0)21} we have $\varphi((-a-b-1,0)) \neq (b-a +1,0 ) $, which is a contradiction. Hence $\varphi((3-a-b,0))\neq (2b+3,2)$. 
\end{itemize} 
Hence $\varphi((2b-a-1,1))\neq (2b-2,0)$.

\textbf{Case C(b)(2)(e):}  If $\varphi((2b-a-1,1))=(2-a-b,2)$ then $\varphi((3-a-b,0)) \sim \varphi((2b-a-1,1))$ and $\varphi((2a-1,2))=(2a-b,1)$ imply
\begin{equation}\label{(3-a-b,0),22}
\varphi((3-a-b,0)) \in \lbrace (2a-b-2,1), (2b-1-a \pm b,0) \rbrace.
\end{equation}
From the Equations \ref{(3-a-b,0),20} and \ref{(3-a-b,0),22} we have
\begin{itemize}
	\item either $2b-1-a \pm b=3b+a \pm 1$, i.e., $2a+2b=0$ or $2a=0$ or $2a+2=0$, i.e., $2a+2b=0$, which are impossible by the Lemma \ref{relations} but $2a+2b+2=0$ may hold. We have this Equation when $\varphi((a+b+1,0))= (b-a +1,0 ) $. So $2a+2b+2=0$ imply $a+b+1=-a-b-1$, hence $\varphi((-a-b-1,0))= (b-a +1,0 ) $. But from the Equations \ref{(-a-b-1,0)20} and \ref{(-a-b-1,0)21} we have $\varphi((-a-b-1,0)) \neq (b-a +1,0 ) $, which is a contradiction. Hence $\varphi((3-a-b,0)) \neq (3b+a+1,0)$.
	\item or $2a-b-2=4a \pm b$, i.e., $2a+2=0$, i.e., $2a+2b=0$ which is impossible by the Lemma \ref{relations} but $2a+2b+2=0$ may hold. We have this Equation when $\varphi((a+b+1,0))= (b-a +1,0 ) $. So $2a+2b+2=0$ imply $a+b+1=-a-b-1$, hence $\varphi((-a-b-1,0))= (b-a +1,0 ) $. But from the Equations \ref{(-a-b-1,0)20} and \ref{(-a-b-1,0)21} we have $\varphi((-a-b-1,0)) \neq (b-a +1,0 ) $, which is a contradiction, so $\varphi((3-a-b,0)) \neq (4a+b,1)$.
\end{itemize}  Hence $\varphi((2b-a-1,1))\neq (2-a-b,2)$.

\textbf{Case C(b)(2)(f):}  If $\varphi((2b-a-1,1))=(2-a+b,2)$ then $\varphi((3-a-b,0)) \sim \varphi((2b-a-1,1))$ and $\varphi((2a-1,2))=(2a-b,1)$ imply
\begin{equation}\label{(3-a-b,0),23}
\varphi((3-a-b,0)) \in \lbrace (2a-b+2,1), (2b-1+a \pm b,0) \rbrace.
\end{equation}
From the Equations \ref{(3-a-b,0),20} and \ref{(3-a-b,0),23} we have,
\begin{itemize}
	\item either $2a-b+2=4a \pm b$, i.e., $2a-2=0$, i.e., $2a-2b=0$ or $2a+2b-2=0$ which are impossible by the Lemma \ref{relations}.
	\item or $2b-1+a \pm b=3b+a \pm 1$, this gives rise to four equations, out of three are impossible by the Lemma \ref{relations}, namely $2=0$, $2b=0$, $2b+2=0$, i.e., $2a+2b=0$. The only possibility which remains is $3b+a-1=3b+a-1$, which is an identity. So assuming this to be the case, we have $\varphi((3-a-b,0))=(3b+a-1,0)$, $\varphi((2b-a-1,1))=(2-a+b,2)$ and $\varphi((3a-b,2))=(3a+1,1)$. 
\end{itemize} 

Now $\varphi((4,0)) \sim \varphi((3b,1))=(3,2) $, $\varphi((2,0))=(2b,0)$, $\varphi((3a+b,2))=(3a-1,1)$ and  $\varphi((3a-b,2))=(3a+1,1)$  imply $\varphi((4,0))=(4b,0)$. $\varphi((2b,0))\sim \varphi((1,2))=(b,1)$, $\varphi((0,0))=(0,0)$, $\varphi((a+1,1))=(a-b,2)$ and $\varphi((a-1,1))=(a+b,2)$ imply $\varphi((2b,0))=(2,0)$.

Proceeding this way we can show that $\varphi((2k,0))=(2kb,0) $, forall $k \in \Bbb Z$. So we have $\varphi((2b,0))=(2a,0)$, where $k=a$, which is a contradiction as we have shown that $\varphi((2b,0))=(2,0)$ and $2a \neq 2$ by Lemma \ref{relations}. Therefore we have $\varphi((a+1,1))\neq (a-b,2)$ and $\varphi((a-1,1))\neq (a+b,2)$ .

\textbf{Case D:} Let $\varphi((a+1,1))=(a+b,2)$, $\varphi((a+b+1,0))= (a+b+1,0 ) $, $\varphi((a+b,2))= (a+1,1)$, $\varphi((a-1,1))=(a-b,2)$, $\varphi((b-a+1,0))= (b-a +1,0 ) $, $\varphi((a-b,2))= (a-1,1)$. Now  $\varphi((a+b,2))= (a+1,1)$, $\varphi((a-b,2))= (a-1,1)$ and $\varphi((2,0)) \sim \varphi((b,1))=(1,2)$ imply $\varphi((2,0))=(2b,0)$.

Consider the cycle $C_1: (2,0) \sim (b,1) \sim (a+b,2) \sim (a+b+1,0) \sim (2b+a+1,1) \sim (2a+1,2) \sim (2,0)$. then $\varphi(C_1): (2b,0) \sim (1,2) \sim (a+1,1) \sim (a+b+1,0) \sim \varphi((2b+a+1,1)) \sim \varphi((2a+1,2)) \sim (2b,0)$. $\varphi((2b+a+1,1)) \sim (a+b+1,0)$, $\varphi((a+1,1))=(a+b,2)$ and $\varphi((a+b,2))=(a+1,1)$ imply $\varphi((2b+a+1,1)) \in \lbrace (a+1+2b,1), (2+b+a,2)  \rbrace$. Now $\varphi((2a+1,2)) \sim (2b,0)$ and $\varphi((b,1))=(1,2)$ imply 
\begin{equation}\label{(2a+1),23}
\varphi((2a+1,2)) \in \lbrace (2a \pm b,1), (3,2) \rbrace.
\end{equation}
\textbf{Case D(a):}  If $\varphi((2b+a+1,1))=(a+1+2b,1)$ then $\varphi((2b+a+1,1)) \sim \varphi((2a+1,2))$ and $\varphi((a+b+1,0))=(a+b+1,0)$ imply 
\begin{equation}\label{(2a+1),24}
\varphi((2a+1,2)) \in \lbrace (1+b+2a \pm b,2), (b+a+3,0) \rbrace.
\end{equation}
then From the Equation \ref{(2a+1),23} and \ref{(2a+1),24} we have, $1+b+2a \pm b=3$, i.e., $2a-2=0$, i.e., $2a-2b=0$ or $2a+2b-2=0$, which are impossible by the Lemma \ref{relations}. Hence $\varphi((2b+a+1,1))\neq (a+1+2b,1)$.

\textbf{Case D(b):}  If $\varphi((2b+a+1,1))= (2+b+a,2)$ then $\varphi((2b+a+1,1)) \sim \varphi((2a+1,2))$ and $\varphi((a+b+1,0))=(a+b+1,0)$ imply 
\begin{equation}\label{(2a+1),25}
\varphi((2a+1,2)) \in \lbrace (2a+1+b \pm 1,1), (3b+a+1,0) \rbrace.
\end{equation}
From the Equation \ref{(2a+1),23} and \ref{(2a+1),24} we have, $2a+1+b \pm 1=2a \pm b$, This gives rise four equations out of which three are impossible by Lemma \ref{relations}, namely $2=0$, $2b=0$, $2b+2=0$, i.e., $2a+2b=0$. The only possiblity which remains is $2a+b=2a+b$, which is an identity. So asssuming this to be the case we have $\varphi((2a+1,2))=(2a+b,1)$, $\varphi((2b+a+1,1))= (2+b+a,2)$.

Now consider the cycle $C'_1: (2,0) \sim (b,1) \sim (a+b,2) \sim (1+a-b,0) \sim (2b+1-a,1) \sim (2a-1,2) \sim (2,0)$. Then $\varphi(C'_1):(2b,0) \sim (1,2) \sim (a+1,1) \sim \varphi((1+a-b,0)) \sim \varphi((2a+1-a,1)) \sim \varphi((2a-1,2))$. $\varphi((2a-1,2)) \sim (2b,0)$, $\varphi((b,1))=(1,2)$ and $\varphi((2a+1,2))=(2a+b,1)$ imply $\varphi((2a-1,2)) \in \lbrace (2a-b,1), (3,2) \rbrace.$ Now  $ \varphi((1+a-b,0)) \sim(a+1,1)$, $\varphi((a+b+1,0))=(a+b+1,0)$ and $\varphi((b,1))=(1,2)$ imply $\varphi((1+a-b,0)) \in \lbrace (b+a-1,0), (1+2b,2) \rbrace$. $\varphi((2b+1-a,1)) \sim \varphi((1+a-b,0))$ and $\varphi((a+b,2))=(a+1,1)$ imply
\begin{equation}\label{(2b+1-a,1)23}
\varphi((2b+1-a,1)) \in \lbrace  (a+1-2b,1), (1+b-a \pm 1,2), (a+3,1), (b+2a \pm b,0) \rbrace.
\end{equation}
\textbf{Case D(b)(1):}  If $\varphi((2a-1,2))=(3,2)$ then $\varphi((2b+1-a,1)) \sim \varphi((2a-1,2))$ and $\varphi((2,0))=(2b,0)$ imply
\begin{equation}\label{(2b+1-a,1)24}
\varphi((2b+1-a,1)) \in \lbrace  (4b,0), (3a \pm 1,1)\rbrace.
\end{equation}
From the Equations \ref{(2b+1-a,1)23} and \ref{(2b+1-a,1)24} we have, 
\begin{itemize}
	\item either $4b=b+2a \pm b$, i.e., $2a-2b=0$ or $2a-4b=0$, which are impossible by the Lemma \ref{relations}.
	\item $3a \pm 1= a+1-2b$, i.e., $2a+2b=0$ or $2a+2b-2=0$, which are impossible by the Lemma \ref{relations}.
	\item $3a \pm 1=a+3$, i.e., $2a-2=0$, i.e., $2a-2b=0$ or $2a=4$, i.e., $4a-2b=0$. Now $2a-2b=0$ s impossible by the Lemma \ref{relations} but $4a-2b=0$ holds only if $n=7,14$.
\end{itemize}
Thus $\varphi((2a-1,2))=(3,2)$ only if $n=7,14$.

\textbf{Case D(b)(2):}  If $\varphi((2a-1,2))=(2a-b,1)$ then $\varphi((2b+1-a,1)) \sim \varphi((2a-1,2))$ and $\varphi((2,0))=(2b,0)$ imply
\begin{equation}\label{(2b+1-a,1)25}
\varphi((2b+1-a,1)) \in \lbrace  (2b-2,0), (2-a \pm b,2)\rbrace.
\end{equation}
From the Equations \ref{(2b+1-a,1)23} and \ref{(2b+1-a,1)25} we have, \begin{itemize}
	\item either $2b-2=b+2a \pm b$, i.e., $2a+2=0$, i.e., $2a-2b+2=0$, which are impossible by the Lemma \ref{relations}.
	\item or $2-a \pm b= 1+b -a \pm 1$, This gives rise four equations out of which three are impossible by Lemma \ref{relations}, namely $2=0$, $2b=0$, $2b-2=0$, i.e., $2a-2b=0$. The only possiblity which remains is $2-a+b=2+b-a$, which is an identity. So asssuming this to be the case we have $\varphi((2b+1-a,0))=(2+b-a,0)$, $\varphi((2a-1,2))=(2a-b,1)$ and $\varphi((1+a-b,0))=(b+a-1,0)$.
\end{itemize}

Now $\varphi((3b,1)) \sim \varphi((2,0))=(2b,0) $, $\varphi((2a+1,2))=(2a+b,1)$ and $\varphi((2a-1,2))=(2a-b,1)$ imply  $\varphi((3b,1)) =(3,2)$.

Consider the cycle $C_2: (2,0) \sim ( 3b,1) \sim (3a+b,2) \sim (3+a+b,0) \sim (2b+a+1,1) \sim (2a+1,2) \sim (2,0)$. Then $\varphi(C_2): (2b,0) \sim (3,2) \sim \varphi((3a+b,2) \sim \varphi((3+a+b,0)) \sim (a+b+2,2) \sim (2a+b,1) \sim (2b,0)$. $\varphi((3+a+b,0)) \sim (a+b+2,2) $, $\varphi((a+b+1,0))=(a+b+1,0)$ and $\varphi((2a+1,2))=(2a+b,1)$ imply $\varphi((3+a+b,0)) \in \lbrace (2a+b+2,1), (3b+a+1,0) \rbrace $. $\varphi((3a+b,2)) \sim (3,2)$  and $\varphi((2,0))=(2b,0)$ imply  
\begin{equation}\label{(3a+b,2)20}
\varphi((3a+b,2)) \in \lbrace (4b,0), (3a \pm 1,1) \rbrace.
\end{equation}
\textbf{Case D(b)(2)(a):}  If $\varphi((3+a+b,0))=(2+b+2a,1)$ then $\varphi((3a+b,2) \sim \varphi((3+a+b,0))$ and  $\varphi((2b+a+1,1))=(a+b+2,2)$ imply 
\begin{equation}\label{(3a+b,2)21}
\varphi((3a+b,2)) \in \lbrace (2a+1+2b \pm 1,0), (3b+a+2,2) \rbrace.
\end{equation}
From \ref{(3a+b,2)20} and \ref{(3a+b,2)21} we have, $2a+1+2b \pm 1= 4b$, i.e., $2a-2b=0$ or $2a-2b+2=0$, which are impossible by Lemma \ref{relations}. Hence If $\varphi((3+a+b,0)) \neq (2+b+2a,1)$.

\textbf{Case D(b)(2)(b):}  If $\varphi((3+a+b,0))=(a+1+3b,0)$ then $\varphi((3a+b,2) \sim \varphi((3+a+b,0))$ and  $\varphi((2b+a+1,1))=(a+b+2,2)$ imply 
\begin{equation}\label{(3a+b,2)22}
\varphi((3a+b,2)) \in \lbrace (1+b+3a \pm b,1), (b+a+4,2) \rbrace.
\end{equation}
From \ref{(3a+b,2)20} and \ref{(3a+b,2)22} we have, $1+b+3a \pm b=3a \pm 1$, this give rise to four equations, out of three are impossible by the Lemma \ref{relations}, namely $2=0$, $2b=0$ and $2b+2=0$, i.e., $2a+2b=0$. The only possibility which remains is $1+3a=3a+1$, which is an identity. So assuming to be the case, we have $\varphi((3a+b,2))=(3a+1,1)$ and $\varphi((3+a+b,0))=(3b+a+1,0)$.

Now consider the cycle $C'_2: (2,0) \sim (3b,1) \sim (3a-b,2) \sim (3-a+b,0) \sim (2b-1+a,1) \sim (2a+1,2)$. Then $\varphi(C'_2): (2b,0) \sim (3,2) \sim \varphi((3a-b,2)) \sim \varphi((3-a+b,0)) \sim \varphi((2b-1+a,1)) \sim (2a+b,1) \sim (2b,0)$. $\varphi((3a-b,2)) \sim(3,2)$, $\varphi((a-b,2))=(a-1,1)$ and $\varphi((2,0))=(2b,0)$ imply $\varphi((3a-b,2)) \in \lbrace (3a-1,1), (4b,0) \rbrace$. $\varphi((2b-1+a,1)) \sim (2a+b,1)$, $\varphi((2b+a+1,1))=(a+b+2,2)$ and $\varphi((2,0))=(2b,0)$ imply $\varphi((2b-1+a,1)) \in \lbrace (2+a-b,2), (2b+2,0) \rbrace$. Now $\varphi((3-a+b,0)) \sim \varphi((3a-b,2))$ and $\varphi((3b,1))=(3,2)$ imply
\begin{equation}\label{(3-a+b,0)20}
\varphi((3-a+b,0)) \in \lbrace (3b-a+1,0), (3-2b,2), (5,2), (4a \pm b,1)\rbrace.
\end{equation}
\textbf{Case D(b)(2)(c):} If $\varphi((2b-1+a,1))=(2b+2,0)$ then $\varphi((3-a+b,0)) \sim \varphi((2b-1+a,1))$, $\varphi((3b,1))=(3,2)$ and $\varphi((2a+1,2))=(2a+b,1)$ imply
\begin{equation}\label{(3-a+b,0)21}
\varphi((3-a+b,0)) \in \lbrace (2a+3b,1), (3+2a,2)\rbrace.
\end{equation}
From  the Equations \ref{(3-a+b,0)20} and \ref{(3-a+b,0)21} we have,
\begin{itemize}
	\item either $2a+3b=4a \pm b$, i.e., $2a-2b=0$ or $2a-4b=0$, which are impossible by lemma \ref{relations}.
	\item or $3+2a=3-2b$, i.e., $2a+2b=0$, which is impossible by Lemma \ref{relations}.
	\item or $3+2a=5$, i.e., $2a-2=0$, i.e., $2a-2b=0$,  which is impossible by Lemma \ref{relations}.
\end{itemize}  
Hence $\varphi((2b-1+a,1))\neq (2b+2,0)$.

\textbf{Case D(b)(2)(d):} If $\varphi((2b-1+a,1))=(2+a-b,2)$ then $\varphi((3-a+b,0)) \sim \varphi((2b-1+a,1))$, $\varphi((b-a+1,0))=(b-a+1,0)$ and $\varphi((2a+1,2))=(2a+b,1)$ imply
\begin{equation}\label{(3-a+b,0)22}
\varphi((3-a+b,0)) \in \lbrace (3b+1-a,0), (2a+b-2,1)\rbrace.
\end{equation}
From  the Equations \ref{(3-a+b,0)20} and \ref{(3-a+b,0)22} we have,  
\begin{itemize}
	\item either $2a+b-2=4a \pm b$, i.e., $2a+2=0$, i.e., $2a+2b=0$ or $2a-2b+2=0$,  which are impossible by lemma \ref{relations}.
	\item $3b+1-a=3b-a \pm 1$, this give rise to two equations, out of one is impossible, namely $2=0$. The only possibility which remains is $3b+1-a=3b-a+1$, which is an identity. So assuming to be the case, we have, $\varphi((3-a+b,0))=(3b-a+1,0)$, $\varphi((3a-b,2))=(3a-1,1)$ and $\varphi((2b+a-1,1))=(2+a-b,2)$.
\end{itemize}

Now $\varphi((4,0)) \sim \varphi((3b,1))=(3,2) $, $\varphi((2,0))=(2b,0)$, $\varphi((3a+b,2))=(3a+1,1)$ and  $\varphi((3a-b,2))=(3a-1,1)$  imply $\varphi((4,0))=(4b,0)$. $\varphi((2b,0))\sim \varphi((1,2))=(b,1)$, $\varphi((0,0))=(0,0)$, $\varphi((a+1,1))=(a+b,2)$ and $\varphi((a-1,1))=(a-b,2)$ imply $\varphi((2b,0))=(2,0)$.

Proceeding this way we can show that $\varphi((2k,0))=(2kb,0) $, forall $k \in \Bbb Z$. So we have $\varphi((2b,0))=(2a,0)$, where $k=a$, which is a contradiction as we have shown that $\varphi((2b,0))=(2,0)$ and $2a \neq 2$ by Lemma \ref{relations}. Therefore from the case \textbf{D(b)(1)} we have $\varphi((a+1,1))=(a+b,2)$ and $\varphi((a-1,1))=(a-b,2)$ only if $n=7,14$. This completes the proof. \qed

Now, the proof of Theorem \ref{main-theorem} follows from Theorem \ref{phi(1,2)=(-1,2)}, Theorem \ref{phi(1,2)=(-b,1)} and Theorem \ref{phi(1,2)=(b,1)}.

\section{Open Issues}
In this paper, we introduced an infinite family of half-transitive Cayley graphs. However, a few issues are still pending and can be topics for further research.
\begin{enumerate}
\item {\bf Full Automorphism Group:} It was shown that $\langle \alpha,\beta,\gamma \rangle$ is a subgroup of the full automorphism group. It remains to be shown (as observed in SageMath) that $G=\langle \alpha,\beta,\gamma \rangle$ for $n\neq 7,14$.
\item {\bf Structural Properties of $\Gamma(n,a)$:}  We have computed the girth for some special values of $n$ and shown that $\Gamma(n,a)$ is Hamiltonian if $n$ is odd. However, the girth and Hamiltonicity for general values of $n$ are still unanswered. Similarly, other structural properties like diameter, domination number are few open issues.
\end{enumerate}

\section*{Acknowledgement}
The second author acknowledge the funding of DST-SERB-SRG Sanction no. $SRG/2019/$ $000475$, Govt. of India.


\section*{Appendix: Sage Code for $\Gamma(n,a)$ for $n=7, a=2$}

\begin{verbatim}
n=7
a=2
b=mod(a^2,n)
A=list(var('A_%d' % i) for i in range(n))
B=list(var('B_%d' % i) for i in range(3))
C = cartesian_product([A, B])
V=C.list()
E=[]
Gamma=Graph()
Gamma.add_vertices(V)
for i in range(n):
  for j in range(3):
    E.append(((A[i],B[j]),(A[mod(a*i+1,n)],B[mod(j-1,3)])))
    E.append(((A[i],B[j]),(A[mod(a*i-1,n)],B[mod(j-1,3)])))
    E.append(((A[i],B[j]),(A[mod(b*i+b,n)],B[mod(j+1,3)])))
    E.append(((A[i],B[j]),(A[mod(b*i-b,n)],B[mod(j+1,3)])))
Gamma.add_edges(E)
G=Gamma.automorphism_group()
for f in G:
   if f((A[0],B[0]))==(A[0],B[0]) and f((A[b],B[1]))==(A[1],B[2]) and
   f((A[1],B[2]))==(A[n-1],B[2]):
   print "sucess"
\end{verbatim}

\end{document}